\documentclass[12pt,a4paper]{amsart}

\usepackage{latexsym} 
 
\usepackage[dvips]{graphics}
\usepackage{epsfig}

\usepackage{amssymb}
\usepackage{amsthm}
\usepackage{amsfonts}
\usepackage{amsmath}
\usepackage{amstext}
\usepackage{epic}
\usepackage{eepic}

\numberwithin{equation}{section}

\theoremstyle{plain}%default
\newtheorem{thm}{Theorem}[section] 
\newtheorem{prop}[thm]{Proposition}
\newtheorem{cor}[thm]{Corollary}
\newtheorem{lem}[thm]{Lemma}

\newtheorem{conj}{Conjecture}
\newtheorem{theorem*}{Theorem}[]

\theoremstyle{definition}
\newtheorem{defn}[thm]{Definition}

\newtheorem{example}[thm]{Example}

\newtheorem{ques}{Question}

\theoremstyle{remark}
\newtheorem{rem}[thm]{Remark}

\newcommand{\C}{\mathbb{C}}
\newcommand{\N}{\mathbb{N}}
\newcommand{\R}{\mathbb{R}}
\newcommand{\K}{\mathbb{K}}

\newcommand{\Z}{\mathbb{Z}}

\newcommand{\proj}{\mathbb{P}}

\newcommand{\alphaanal}{\mathcal L }
\newcommand{\alphan}{\mathcal L_n }
\newcommand{\alphai}{\mathcal L_i }
\newcommand{\ix}{\mathcal X }
\newcommand{\zed}{\mathcal Z }
\newcommand{\va}{\mathbf a }
\newcommand{\ixn}{\mathcal X_n }
\newcommand{\ixnp}{\mathcal X_{n,+} }
\newcommand{\ixnm}{\mathcal X_{n,-} }
\newcommand{\bixnp}{\bar {\mathcal X}_{n,+} }
\newcommand{\M}{\mathcal M }
\newcommand{\zetaa}{Z }
\newcommand{\zetap}{Z_+  }
\newcommand{\zetam}{Z_- }
\newcommand{\zetafa}{Z_f }
\newcommand{\zetafp}{Z_{f,+}  }
\newcommand{\zetafm}{Z_{f,-} }
\newcommand{\mzeta}{\widetilde Z }

\newcommand{\chic}{\chi ^c }
\newcommand{\jac}{jac \, }

\newcommand{\inv}{^{-1}}

\newcommand{\half}{\frac{1} 2}

\DeclareMathOperator{\supp}{supp\,}
\DeclareMathOperator{\mult}{mult\,}
\DeclareMathOperator{\ord}{ord_t\,}

\DeclareMathOperator{\grad}{grad\,}

\DeclareMathOperator{\Jac}{J}

\def\accentclass@{7}
\def\makeacc@#1#2{\def#1{\mathaccent"\accentclass@#2 }}
\makeacc@\cir{017}
\newcommand{\Eo}{{\cir E }}
\title{Motivic-type Invariants of Blow-analytic Equivalence}

\author{Satoshi Koike \& Adam Parusi\'nski}

\address {Department of Mathematics, Hyogo University
of Teacher Education, 942-1 Shimokume, Kato, Yashiro,
Hyogo 673-1494, Japan}

\email {koike@sci.hyogo-u.ac.jp}

\address {Departement de Mathematiques, Universit\'e d'Angers,
   2, bd Lavoisier, 49045 Angers cedex, France}

\email{parus@tonton.univ-angers.fr}

\thanks{This research was originated during the conference New 
Developments in Singularity Theory, Cambridge 2000, 
at the Isaac Newton Institute, 
and was subsequently supported by the University of Angers and Grant-in-Aid 
for Scientific Research (No. 13640070) of Ministry of Education, 
Science and Culture of Japan.  We would like to thank these 
institutions for their support and hospitality.}
\subjclass{14B05, 32S15, 57R45}
\begin{document}

\begin{abstract} 
To a given analytic function germ 
$f:(\mathbb{R}^d,0) \to 
(\mathbb{R},0)$, we associate zeta functions $Z_{f,+}$, $Z_{f,-} 
\in \mathbb{Z} [[T]]$, defined analogously to the motivic zeta 
functions of Denef and Loeser.  We show that our zeta functions are 
rational and that they are invariants of the blow-analytic 
equivalence in the sense of Kuo.  Then we use them together with the 
Fukui invariant to classify the blow-analytic equivalence classes of 
Brieskorn polynomials of two variables. Except special series of 
singularities our method classifies as well the blow-analytic 
equivalence classes of Brieskorn polynomials of three variables.  

\bigskip
\noindent
{\bf Resum\'e.} Soit $f:(\mathbb{R}^d,0) \to 
(\mathbb{R},0)$ un germe de fonctions analytiques.  On associe \`a $f$
des fonctions zeta $Z_{f,+}$, $Z_{f,-} 
\in \mathbb{Z} [[T]]$ d\'efinies de mani\`ere similaire que les fonctions 
zeta motiviques de  Denef et Loeser.  On montre que ces fonctions sont
rationnelles et ne dependent que de la classe d'\'equivalence
blow-analytique au
sens de Kuo de $f$.  En utilisant ces fonctions zeta et l'invariant de
Fukui 
on donne une classification 
des polyn\^omes de Brieskorn de deux variables 
\`a \'equivalence blow-analytique
pr\`es.  Pour les polyn\^omes de Brieskorn de trois variables on
obtient une classification presque compl\`ete.   

\end{abstract}

\maketitle

%%%%%%%%%%%%%%%%%%%%%%%%%%%%%%%%%%%%%%%%%%%%%%%%%%%%%%%%%%%%%%%%%%%%

\vspace{0.4 truecm}  

%{\it Acknowledgment.} 

%%%%%%%%%%%%%%%%%%%%%%%%%%%%%%%%%%%%%%%%%%%%%%%%%%%%%%%%

In this paper we develop techniques that allow us to study and 
distinguish different blow-analytic classes of analytic function
germs $f:(\mathbb{R}^d,0) \to (\mathbb{R},0)$.   For this we
adapt and apply to the real analytic set-up the ideas coming
from  motivic integration, in particular the concept of 
motivic zeta function due to Denef and Loeser.

The notion of blow-analytic equivalence was introduced
by T.-C. Kuo \cite{kuo1} and \cite{kuo2}.
Recall briefly that analytic function germs 
$f, g : (\R^d,0) \to (\R,0)$ are {\em blow-analytically equivalent} 
if there exist real modifications
$\mu : (M,\mu^{-1}(0)) \to (\R^d,0)$,
$\mu^{\prime} : (M^{\prime},\mu^{\prime -1}(0)) \to (\R^d,0)$
and an analytic isomorphism $\Phi : (M,\mu^{-1}(0)) \to
(M^{\prime},\mu^{\prime -1}(0))$ which induces a homeomorphism 
$\phi : (\R^d,0) \to (\R^d,0)$ such that $f = g \circ \phi$.
In this paper we suppose additionally that $\mu$, resp. $\mu'$, is
an isomorphism over the complement of $f\inv (0)$, resp. $g\inv
(0)$.  The blow-analytic equivalence is interesting because 
it does not allow continuous moduli for families of isolated
singularities cf. \cite{kuo2}, and it preserves a deep information 
on the algebraic
structure of the singularity.  For  real singularities,
unlike for the complex ones, the topological
classification is too crude, e.g. $x^{2k} + y^{2m}$ and 
$x^{2n} + y^{2l}$ are always topologically equivalent.  
The blow-analytic equivalence was invented to overcome this 
problem.  
Moreover, as follows from various examples, the blow-analytic
equivalence of real analytic function germs behaves in a similar
way to 
the topological equivalence in the complex case, though there is
no precise result in this direction.  This observation seems to
be confirmed by the main results of this paper.

There exist various criteria of
blow-analytic triviality of families of analytic function germs,
based mainly on toric equi-resolutions 
\cite{fukuiyoshinaga}, \cite{fukuipaunescu}, \cite{abderrahmane},
but there were till now very few results allowing to
distinguish different blow-analytic types and hence to attempt a
classification even in the simplest cases.
The only known up to now invariant of blow-analytic equivalence 
was introduced by Fukui in
\cite {fukui},  see also section \ref{Fukui} below.
In this paper we introduce new invariants that allow us to start 
such a classification.  

\smallskip
The main results of this paper are the following.  
In section \ref{MZF} we associate
to each real analytic function germ $f:(\mathbb{R}^d,0) \to 
(\mathbb{R},0)$ its zeta functions:  $\zetafa , \zetafp , \zetafm  
\in \mathbb{Z} [[T]]$.   We show that they are blow-analytic
invariants in section \ref{Zetainvariant}.  In order to compute 
the zeta functions we propose the formulae in terms of a
resolution (Denef\&Loeser formulae), see section \ref{MZF}, and 
the Thom-Sebastiani Formulae in section \ref{TS}.  
Sections \ref{2BP} and \ref{3V} contain classification results, in
particular a complete classification of blow-analytic types  of 
Brieskorn polynomials of two variables and  a partial
classification in three dimensional case.  

%\smallskip
Our main idea of construction of new invariants 
is based on the following simple observations.
Suppose that $f, g : (\R^d,0) \to (\R,0)$ are blow-analytically equivalent
via a (blow-analytic) homeomorphism $\phi : (\R^d,0) \to
(\R^d,0)$, $f = g \circ \phi$.  Then, firstly, $f$ and $g$ admit 
isomorphic resolutions.  Secondly,  let $\alphaanal (\R^d,0)$
denote the set of germs of analytic arcs at the origin
in $\R^d$.  Then $\varphi$ induces a bijection $\varphi_* : \alphaanal
(\R^d,0) \to \alphaanal (\R^d,0)$ by composition
$\varphi_*(\gamma(t)) = (\varphi \circ \gamma)(t)$.  
In section \ref{MZF} below, using the integration with respect to
the Euler characteristic with compact supports 
on these sets of arcs, we associate
to each real analytic function germ $f:(\mathbb{R}^d,0) \to 
(\mathbb{R},0)$ its zeta functions:  $\zetafa , \zetafp , \zetafm  
\in \mathbb{Z} [[T]]$.   Here we follow the path introduced by 
Denef and Loeser \cite{denefloeser-1},
\cite{denefloeser3}, and inspired by work of Kontsevich \cite{kontsevich}.
The zeta function of Denef and Loeser, and the related topological zeta 
function  cf. \cite{denefloeser0},
provides an important information on the local topology of
complex analytic function germs, see a new proof of 
Thom-Sebastiani
theorem for the Hodge spectrum \cite{denefloeser2} or works on the monodromy 
conjecture, see for instance
\cite{denefloeser3}, \cite {veys}.  
We refer the reader to the survey \cite{denefloeser0} for more information
on the Denef and Loeser construction and its applications.

In section \ref{Zetainvariant} we show that our  zeta functions
are invariants of blow-analytic equivalence in the sense of Kuo.  
The proof is based on formulae \eqref{DL}, \eqref{DLpm}, analogous to 
the formulae of Denef and Loeser, that express the zeta functions of $f$ 
in terms of a resolution.  These formulae are proven by a version of
Kontsevich's change of variable formula, Corollary \ref{change}.  Note
that these results do not follow automatically from the analogous 
ones in the algebraic case, due to the necessity of working with 
non-compact subanalytic sets.  This difficulty is overcomed 
thanks to  the {\L}ojasiewicz's theory of relatively 
semi-algebraic, semi-analytic sets \cite{lojasiewicz}.

Thom-Sebastiani
Formulae, showed in  section \ref{TS}, express the zeta functions of
$f(x) + g(y)$ in terms of the ones of $f$ and $g$.  They have 
interesting consequences.  For instance we get a suspension
property: if the zeta functions of $x^m + g_1(y)$ and $x^m +
g_2(y)$, $m$ even,  coincide then so do the zeta functions of 
$g_1$ and $g_2$.  One may speculate 
that if  $x^2 + g_1(y)$ and $x^2 + g_2(y)$ are blow-analytically
equivalent so are $g_1$ and $g_2$ (this is for instance the
case if we know that the zeta functions distinguish the blow-analytic
types of $g_1$ and $g_2$,
we use this in some special cases).  We do not know the answer to
this question.  We use the Thom-Sebastiani
Formulae
to compute the zeta functions for all Brieskorn
polynomials $f(x_1,\ldots,x_d)= \pm x_1^{p_1} \pm \cdots \pm x_d^{p_d}$.
In section \ref{2BP} we compute the blow-analytic
equivalence classes of 
Brieskorn polynomials of two variables and in section \ref{3V}
most of the  equivalence classes of Brieskorn polynomials of
three variables.  This classification differs from the analytic one. 
For instance, thanks to a  phenomenon typical for real
algebraic geometry, the functions
$x^p +y^{kp}$ and $x^p -y^{kp}$, $p$ odd, $k$ even, are blow-analytically 
equivalent
but not analytically equivalent (over real numbers).  

As we mentioned before the blow-analytic equivalence behaves in a
similar way to the topological equivalence of complex analytic
function germs.  Consider for instance the following example.  
The germs at the origin of $f(x,y,z)= x^3+xy^5 +z^3$ and 
$g(x,y,z)= x^3 + y^7 +z^3$ are not topologically equivalent as 
complex germs.  One may show that any complex analytic function germ 
with the $6$th jet equal to $f$ is
topologically equivalent either to $f$ or $g$, thus there are
exactly two possible topological types.  On the other hand any
real analytic function germ with the $6$th jet equal to $f$ is 
blow-analytically equivalent either to $f$ or $g$.  Of course, 
$f$ and $g$ as real analytic functions germs are topologically
equivalent (they are equivalent to a nonsingular germ).  We show
in subsection \ref{nB} that $f$ and $g$ are not blow-analytically
equivalent.

%\smallskip
Due to the presence of some phenomena typical for the real algebraic 
geometry it is interesting to compare the properties of our zeta functions
to the ones of Denef and Loeser.  
For instance our sign zeta functions, $\zetap, \zetam$, 
correspond to the monodromic zeta function of Denef and Loeser, a
phenomenon similar to the one studied in \cite{mccroryparusinski} in 
a different context. Note also 
that our zeta functions are not really motivic and have only
integer coefficients.  This is due to the fact that
the Euler characteristic with compact
support is the only numerical invariant of the semi-algebraic motifs 
as defined topologically in \cite{quarez}.  

Moreover the zeta functions introduced in this paper do not
distinguish all classes of blow-analytic equivalence and we
are far from a complete classification even in the weighted
homogeneous non-degenerate case.  This problem may  be attack by
hunting new motivic  invariants in the real algebraic, and not 
semi-algebraic, 
set-up.  Even if one knows such
invariants it is not clear whether one can apply them to study
the equivalence that is merely blow-analytic (and not
``blow-algebraic").  On the other hand there is a variety of work
in real algebraic and analytic geometry related to the space of 
analytic arcs that can be probably approached by the techniques
of motivic integration, cf.  \cite{kurdyka1}, \cite{bierstonemilman0}, 
\cite{kurdyka2}.

%\smallskip
We finish the introduction with more precise questions.  
Let $f$, $g : (\C^n,0) \to (\C,0)$ be weighted homogeneous
polynomials with isolated singularities.
It is known after  \cite{nishimura},
\cite{saeki}, \cite{yau}, \cite{yoshinagasuzuki}, for $n = 2$, $3$, 
that if  $(\C^n,f^{-1}(0))$ and $(\C^n,g^{-1}(0))$
are homeomorphic as germs at $0 \in \C^n$,
then their systems of weights coincide.
We propose the following corresponding question.
\begin{ques}\label{baweights}
Let $f$, $g : (\R^n,0) \to (\R,0)$ be weighted homogeneous
polynomials with isolated singularities.
Suppose that $f$ and $g$ are blow-analytically equivalent.
Then, do their systems of weights coincide?
\end{ques}

Let $\K = \R$ or $\C$,
and let $J^r_{\K}(n,1)$ denote the set of $r$-jets of
analytic function germs $(\K^n,0) \to (\K,0)$.
We identify $r$-jets with polynomial representatives 
of degree not exceeding $r$.
We say that $w \in J^r_{\K}(n,1)$ satisfies {\em the Kuiper-Kuo
condition} (\cite{kuiper}, \cite{kuo0}) if there are 
$C$, $\alpha > 0$, such that
$$
|\grad w(x)| \ge C|x|^{r-1} \ \ \text { for } \ \ |x| < \alpha.
$$
Concerning blow-analytic sufficiency of jets, T.-C. Kuo
gave the following conjecture and has affirmatively 
proved it in the two variables case.

\begin{conj}\label{basufficiency}
Let $w \in J^r_{\R}(n,1)$.
Suppose that $w$ satisfies the Kuiper-Kuo condition
as a complex $r$-jet.
Then $w$ is blow-analytically sufficient in 
$C^{\omega}$-functions.
\end{conj}

\smallskip
\noindent \emph{Convention:}
By the Brieskorn polynomials of $d$ variables we mean 
$f(x_1,\ldots, x_d) = a_1x_1^{p_1} + a_2x_2^{p_2} +\cdots
+a_dx_d^{p_d}$, 
$a_i\ne 0$.  
Since their analytic types depend only on the signs of $a_i$, in 
order to simplify the notation, in this paper we consider only the 
Brieskorn polynomials of the form  $f(x_1,\ldots, x_d) = 
\pm x_1^{p_1} \pm x_2^{p_2} \pm  \cdots \pm x_d^{p_d}$.  
\medskip

%%%%%%%%%%%%%%%%%%%%%%%%%%%%%%%%%%%%%%%%%%%%%%%%%%%%%%%%%%%%%%
\bigskip
\section{Motivic zeta function of analytic function germ}
\label{MZF}
\medskip

\subsection{Definition of the zeta functions}
Consider the space of analytic arcs at the origin $0\in 
\R^d$   
$$
\alphaanal = \alphaanal (\R^d,0) := \{ \gamma: (\R,0)\to (\R^d,0);
\gamma \text { analytic} \} 
$$
and the one of truncated arcs 
 $$
\alphan := \{ \gamma \in \alphaanal; \gamma(t) = 
\va _1 t + \va _2 t^2 + \cdots +  \va _n t^n , \va_i\in \R^d \}.  
$$
Given an analytic function $f:(\R^d,0)\to (\R,0)$.
For $n\ge 1$ we denote
\begin{eqnarray*} 
& \ixnp(f) := & \{ \gamma\in \alphan ; f\circ \gamma =
ct^n + \cdots ,c>0 \},  \\
& \ixnm (f):= & \{ \gamma\in \alphan ; f\circ \gamma =
ct^n + \cdots , c< 0\} , \\
& \ixn (f) := &  \{ \gamma\in \alphan ; f\circ \gamma =
ct^n + \cdots , c\neq 0\} . 
\end{eqnarray*}
We define the { \it positive, negative, and total zeta function} 
of $f$ by 
\begin{eqnarray*} 
& \zetafp (T) := & \sum_{n\ge 1}  \,(-1) ^{-nd} \chic(\ixnp )  T^n,  
\\
& \zetafm (T):= & \sum_{n\ge 1}  \,(-1) ^{-nd} \chic(\ixnm )  T^n, \\
& \zetafa (T):= & \sum_{n\ge 1}  \, (-1) ^{-nd} \chic(\ixn )  T^n = 
\zetafp (T)+ \zetafm (T) , 
\end{eqnarray*} 
where $\chic$ denotes the Euler characteristic with compact supports.  
If $f$ is fixed we shall often drop $f$ and write simply
$\ixnp $ for $\ixnp (f)$, $\zetap$ for  $\zetafp$,
and so on. 

\medskip
\begin{rem}\label{trivial1}
The map $\varphi: \ixn \to \R^*$ that associates to
$\gamma$ the first non-zero coefficient of $f\circ \gamma$,
that is $\varphi 
(\gamma) = c$ if $ f\circ \gamma = ct^n + \cdots$,
 is a trivial fibration over  $\R_{<0}$ and $\R_{>0}$ 
(for $n$ odd it is trivial over $\R^*$).
This can be easily shown
using the following action of $\R^*$  
$$ 
\varphi (\gamma (\alpha t)) = \alpha ^n \varphi (\gamma(t)) ,
\quad  \alpha \in \R^* . 
$$
\end{rem}

\medskip
\begin{rem}\label{motifs}
Our zeta function is an incarnation of the motivic zeta function of
Denef \& Loeser \cite{denefloeser0}, \cite{denefloeser2},
\cite{denefloeser3}.  Instead of using the algebraic motifs we use
just the Euler characteristic with compact supports that is the Euler
characteristic of the sheaf cohomology with compact supports, with
coefficients in the constant sheaf $\Z$.   By the long exact
cohomology sequence of the pair it satisfies the
following additivity property: for all 
locally compact semialgebraic $A$ and $B$, $B$ closed in $A$, $\chic
(A) = \chic (A\setminus B) + \chic (B)$.  One may show that $\chic$ is
the only topological invariant of semi-algebraic sets additive in this
sense, cf.  \cite{quarez}. 
\end{rem}

\bigskip
\subsection{Denef \& Loeser's formulae}\label{DLformulae}

Let $\sigma: (M,\sigma \inv (0))\to (\R^d,0)$ 
be a modification of $\R^d$
such that $f\circ \sigma$ and the jacobian determinant 
$\jac \sigma$ of $\sigma$ are normal crossings
simultaneously (we may define $\jac \sigma$ locally
using any local system of
coordinates on $M$).  For instance if $\sigma$ is a composition
of blowings-up with smooth centers that are in normal crossings
with the old exceptional divisors then $\jac \sigma$
is normal crossings.  We also assume that 
$\sigma$ is an isomorphism over the complement of the 
zero set of $f$.    The existence of such a modification
is guaranteed by \cite{hironaka}, \cite{bierstonemilman}.  
We denote by $E_i$, $i\in J$, the irreducible components of 
$(f\circ \sigma) \inv (0)$ (in $\sigma \inv (B_\varepsilon)$,
where $B_\varepsilon$ is a small ball in $\R^d$ centered at
the origin).  We may also suppose
that $\sigma \inv (0)$ is the union of
some of $E_i$.   
For each $i\in J$ we denote $N_i= \mult _{E_i} f\circ \sigma$ and 
$\nu_i= \mult _{E_i} \jac \sigma +1$.  
Denote for $i\in I$ and $I\subset J$, 
$\Eo _i = E_i \setminus \bigcup_{j\neq i} E_j$,  
$E_I = \bigcap _{i\in I} E_i$, $\Eo _I = E_I
\setminus \bigcup _{j\in  J\setminus I} E_j$.
Using the Kontsevich formula of change of 
variables in the motivic integral \cite{kontsevich},
\cite{denefloeser1}, \cite{looijenga} we shall show in section 
\ref{Zetainvariant} that 
\begin{equation}\label{DL}
Z(T) = \sum_{I\neq \emptyset} (-2)^{|I|} \chic (\Eo _I \cap 
\sigma\inv (0)) \prod_{i\in I} \frac {(-1)^{\nu_i} T^{N_i}} 
{1- (-1)^{\nu_i} T^{N_i}} .  
\end{equation} 
Let $\Eo _{I,k}$ be a connected component of $\Eo_I$ and let 
$x\in \Eo _{I,k}$.  Then, near $x$, the complement of $(f\circ
\sigma) \inv (0)$ consists of $2^{|I|}$ chambers, $f$ being non-zero 
on each of them.  Denote by $\alpha_+ (\Eo _{I,k})$ , resp. 
$\alpha_- (\Eo _{I,k})$, the number of such chambers where $f\circ
\sigma$ is positive, resp. negative.
Again using  Kontsevich's argument one gets   
\begin{equation}\label{DLpm}
Z_\pm (T) = \sum_{I\neq \emptyset} (-1)^{|I|} 
 \Bigl (\sum_k  
\alpha_\pm (\Eo _{I,k}) \chic (\Eo _{I,k} \cap  \sigma
\inv (0)) \Bigr )
\prod_{i\in I} \frac  {(-1)^{\nu_i}
T^{N_i}} {1- (-1)^{\nu_i} T^{N_i}}  .  
\end{equation} 

The formulae \eqref{DL}, \eqref{DLpm} will be shown in section 
\ref{Zetainvariant} below.

\bigskip
\subsection{Examples} \label{Examples}
\subsubsection{}\label{xtom}
$f(x) = x^m$, $x\in \R$, $m>0$.  Then 
\begin{equation}
\ixn = 
\begin{cases} \{ \gamma (t) = a_k t + \cdots + a_n t^n; 
a_k \ne 0\} \simeq \R^* \times \R^{n-k} \quad \text { if }
n =km  \\
\emptyset \hfill \quad \text { otherwise} .
\end{cases}
\end{equation}
That is $
\chic (\ixn ) = (-2) \chic (\R^{n-k}) = (-2) (-1)^{k(m-1)}
\quad \text { if } n=km,$
and
$$
Z(T) = \sum _{n=km>0} (-1)^{km} (-2) (-1)^{k(m-1)} T^{km}
= 
2(T^m-T^{2m}+T^{3m} - \cdots ) .  
$$
Of course, the same formula can be obtained by \eqref{DL} by
 taking $\sigma$ equal to the identity
$$
Z(T) = (-2) \frac {-T^m}{1+T^m} .
$$
If $m$ is odd then $Z_+(T) = Z_-(T) = \half Z(T)$.
If $m$ is even then $Z_+(T) = Z(T)$, $Z_-(T) = 0$.

\subsubsection{}\label{2k2k}
$f(x,y) = x^{2k} + y^{2k}$, $(x,y) \in \R^2$.
We may desingularize $f$ by blowing-up the origin with the
 exceptional divisor $\proj ^1$.
Since $\chic (\proj^1) =0$
 we get by \eqref{DL} and \eqref{DLpm}  
$$
 Z_+(T) = Z_-(T) = Z(T)=0.
$$

\subsubsection{}\label{2-2}
$f(x,y) = x^{2} - y^{2}$, $(x,y) \in \R^2$.  Since
$f$ is already normal crossing we apply \eqref{DL} to $\sigma = id$.  Then
$$
Z(T) = (-2)^2 \chic (point) \frac {-T}{1+T} \frac {-T}{1+T}
= 4\frac {T^2}{(1+T)^2} = 4T^2(1-2T+3T^{2} - \cdots ).
$$

\subsubsection{}\label{mm}
$f(x,y) = x^{m} + y^{m}$, $(x,y) \in \R^2$, $m$ odd.
Then $f$ can be desingularized by one blowing-up with the
exceptional divisor $\proj^1$.  Now,
$(f\circ \sigma)\inv (0)$ contains as well the strict
transform of $f\inv (0)$ that is a smooth curve meeting
the exceptional divisor transversally at a point.
Hence
\begin{eqnarray*}
& & Z(T)  =  (-2) (-1)
\frac {T^m}{1 - T^m} + (-2)^2 
\frac {T^m}{1 - T^m} \frac {-T}{1+T} \\
& & =
2T^m(1-2T +2T^2 - \cdots + 2T^{m-1}- T^{m} + 0 + \cdots
+ 0 + T^{2m} - 2 T^{2m+1} + \cdots ).  
\end{eqnarray*}
Clearly $Z_+(T)= Z_-(T) = \half Z(T) $.

%\bigskip
\subsubsection{Zeta functions of a product.}
\label{Product} 

Let $f(x,y):(\R^d,0) \to (\R,0)$,  $f(x,y) = f_1(x)f_2(y)$ where
$f_i: (\R^{d_i},0) \to (\R,0)$, $i=1,2$. 
Then it is easy to check the following formulae. 
\begin{equation*}
Z_f = Z_{f_1} Z_{f_2}, \quad Z_{f,+} = Z_{f_1,+} Z_{f_2,+} + 
Z_{f_1,-} Z_{f_2,-},\quad  Z_{f,-} = Z_{f_1,+} Z_{f_2,-} + 
Z_{f_1,+} Z_{f_2,-}.
\end{equation*}
If $f(x,y) = f_1(x)f_2(y)$, $f_1(0) =0$ but $f_2(0)> 0$, then 
$Z_{f,\pm} = Z_{f_1,\pm}$ and the signs are swapped if $f_2(0)<0$.  
  
Let $f(x)= u(x)\prod _{i=1}^k x_i^{N_i}$, $N_i\ge 1, k\ge 1$,
with $u(0) \neq 0$.  
By the above and  
example \ref{xtom} 
\begin{equation}\label{ncrossing}
\zetafa (T) = (-2)^k \prod _{i=1}^k \frac {-T^{N_i}} {1+T^{N_i}}. 
\end{equation}
If one of $N_i$ is odd then $\zetafp (T) = \zetafm (T)$.
If they are all even and $u(0) <0$, resp. $u(0) >0$, then
$\zetafp (T) \equiv 0$, resp. $\zetafm (T)\equiv 0$.

%%%%%%%%%%%%%%%%%%%%%%%%%%%%%%%%%%%%%%%%%%%%%%%%%%%%%%%%%%%%

\bigskip
\section{Thom-Sebastiani formulae}
\label{TS}
\medskip

The Thom-Sebastiani Formulae  
express the zeta functions of $f(x) + g(y)$ in terms
of the zeta functions of $f(x)$ and $g(y)$, $x\in \R^{d_1},
y\in \R^{d_2}$.  
We denote
$(f*g)(x,y) := f(x) + g(y)$.  Fo motivic zeta functions similar
formulae were proposed in 
\cite{denefloeser2}.  

In what follows we denote  
$$
Z_{f,\pm } (T) = \sum a_i^\pm T^i, \quad 
Z_{g,\pm } (T) = \sum b_i^\pm T^i, \quad 
Z_{f*g,\pm } (T) = \sum c_i^\pm T^i. 
$$
Then
$$
Z_{f} (T) = \sum a_i T^i, \quad 
Z_{g} (T) = \sum b_i T^i, \quad 
Z_{f*g} (T) = \sum c_i T^i,
$$
where  $a_i = a_i^+ + a_i^-$ and so on.
Let $A_n = 1 - \sum_1^n a_i, n \ge 1$, $A_0=1$.  Then
$\sum_{i\ge 0} A_i T^i = \frac {1-Z_f(T)}{1-T}$.
Similarly we define $B_n$, $n \ge 0$.

\smallskip
\begin{thm}\label{Thom-Seb}
\begin{eqnarray}\label{c+}
& c_n^+ & = a_n^+b_n^+ + a_n^+ B_n + A_n b_n^+ +
\sum_1^n (-1) ^{n-i} ( a_i^+b_i^- +  a_i^- b_i^+), \\
\label{c-}
& c_n^- & = a_n^-b_n^- + a_n^- B_n + A_n b_n^- +
\sum_1^n (-1) ^{n-i} ( a_i^+b_i^- +  a_i^- b_i^+) ,\\
\label{c} 
& c_n & = a_n^+b_n^+ + a_n^-b_n^- + a_n B_n + A_n b_n +
2 \sum_1^n (-1) ^{n-i} ( a_i^+b_i^- +  a_i^- b_i^+). 
\end{eqnarray}
\end{thm}

Note that, in general, the total zeta function $Z_{f*g} (T)$
depends on all, that is also on the positive and negative zeta
functions of $f$ and $g$ and not only on $Z_{f} (T)$  and
$Z_{g} (T)$ as the following example shows.

\begin{example}\label{example1}
Let $f(x)=x^2$, $g (y) = y^2$.  The zeta functions
of $f$ and $g$ are computed in Subsection \ref{DLformulae}.
The coefficients $A_i$ are given by 
$$
\sum A_i T^i = \frac {1+T}{1+T^2} = 1+T -T^2 -T^3 +
T^4 +T^5 - \cdots.
$$
One may compute easily the zeta functions of $f*g$ using 
theorem \ref{Thom-Seb}.  They, of course, coincide with the
ones given by \ref{2k2k}.  The total zeta function of  
$h(y) = -y^2$ equals that of $g$.  But 
the total zeta functions $f*g$ and $f*h$ are
different, see \ref{2-2}.  
\end{example}

In general formulae \eqref{c+}-\eqref{c} are not easy
to use.  Moreover they are term by term formulae.  If the
zeta functions of $f$ and $g$ are
given by rational functions of $T$, then theorem
\ref{Thom-Seb} does not give a similar form for the zeta
functions of $f*g$.  The Thom-Sebastiani Formulae can be
simplified considerably by introducing {\it the modified zeta
functions} given by   
$$
\mzeta_{f,+} (T) = \sum_{n\ge 1} \tilde{A}^+_n T^n, \quad 
\mzeta_{f,-} (T) = \sum_{n\ge 1} \tilde{A}^-_n T^n, 
$$
where $\tilde{A}^+_n = A_n + a^+_n, \ 
\tilde{A}^-_n = A_n + a^-_n$. Then
\begin{equation}\label{zetatotilde}
\mzeta_{\pm}(T) = \frac {1 - Z(T)}{1 - T} - 1 + Z_{\pm} (T) 
\end{equation}
and if we introduce {\it the total modified zeta function} by 
$\mzeta (T) := \mzeta_{-}(T)  + \mzeta_{+}(T)$ then  
$$
\frac {1-Z(T)}{1-T} = \frac {1+ \mzeta(T)}{1+T}.
$$
We can compute the zeta functions from the modified ones
by the inverse formula
\begin{equation}\label{mzetatozeta}
{Z}_{\pm}(T) = \frac {1 + \mzeta (T)}{1 + T} + 1 + \mzeta_{\pm} (T) .
\end{equation}  

Let
\begin{equation}\label{mzeta}
\mzeta_{f,\pm}(T) = \sum_{i \ge 1} \tilde{A}_i^{\pm} T^i, \, 
\mzeta_{g,\pm}(T) = \sum_{i \ge 1} \tilde{B}_i^{\pm} T^i, \, 
\mzeta_{f*g,\pm}(T) = \sum_{i \ge 1} \tilde{C}_i^{\pm} T^i
\end{equation}
(same signs). The following formulae are equivalent to those of
theorem \ref{Thom-Seb}.

\medskip
\begin{thm}\label{Thom-Seb2}
\begin{equation}\label{tildeC}
\tilde{C}_n^+ = \tilde{A}_n^+ \ \tilde{B}_n^+,
\ \ \tilde{C}_n^- = \tilde{A}_n^- \ \tilde{B}_n^-.
\end{equation}
\end{thm}

\medskip
\begin{example}\label{example2} \hfil

\begin{enumerate}
\item [(a)]
Let $f(x)=x^m$, $m$ odd.  Then,
$$
\mzeta_{f,+}(T)= \mzeta_{f,-}(T)=
T + T^2 + \cdots + T^{m-1} -T^{m+1} -\cdots -
T^{2m-1} + T^{2m+1} +\cdots .
$$
In particular, $\tilde A ^+_{n} = \tilde A ^-_{n} =0$
for $n\in m\N$.

\item [(b)]
Let $f(x)=x^m$, $m$ even.  Then,
\begin{eqnarray}\label{xtomeven}
& \mzeta_{f,+}(T) & = 
T + T^2 + \cdots + T^{m} -T^{m+1} -\cdots -
T^{2m} + T^{2m+1} +\cdots \\
& \mzeta_{f,-}(T)& =
T + T^2 + \cdots + T^{m-1} - T^{m} -\cdots -
T^{2m-1} + T^{2m} +\cdots. 
\end{eqnarray}  
\end{enumerate}
\end{example} 

\begin{cor} \label{xmeven}
Let $f(x) = x^{m}$ or $-x^m$, $m$ even.
Then $Z_{g ,\pm}(T)$ can be computed
from $Z_{f*g ,\pm}(T)$.
\end{cor}

\begin{proof}
As follows from theorem \ref{Thom-Seb2} this suspension
 property holds for any function $f(x)$ for
which all $\tilde A^\pm _{n}$ are non-zero.  This holds for
$f(x) =  \pm x^m$, $m$ even, by example \ref{example2} (b).  
\end{proof}

If $f(x) = x^{m}$, $m$ odd, then, in general, $Z_{g }(T)$
cannot be computed from $Z_{f*g }(T)$.  Nevertheless we have the
following result.

\begin{prop} \label{xmodd}
Let $f(x) = x^{m}$, $m>1$ odd, and let $g(y) = \pm y^{k}$.  Then 
$k$ is determined by the zeta functions of $f*g$.  If, moreover, $k$ 
is even and not divisible by $m$ then the sign at $y^{k}$ 
is determined by the zeta functions of $f*g$.
\end{prop}

\begin{proof}
We use notation \eqref{mzeta} for the modified zeta functions of 
$f$, $g$, and $f*g$.   
If $\tilde C^+_n =0$ for  $n\notin m\N$  
then, by theorem \ref{Thom-Seb2}, $\tilde B^+_n =0$.  Then $k$ is 
odd and equals the minimum of such $n$.  Similarly, if
there is $n\notin m\N$ such that
$\tilde C^+_n \ne \tilde C^-_n$
then,  $k$ is even and equals the minimum of such $n$. 
Thus suppose that 
$$
\tilde B^+_n = \tilde B^-_n \ne 0 \qquad \text{for all } n\notin m\N. 
$$
Then $k$ is a multiple of $m$ and equals the minimal $n=pm$ that
produce a sign change $\tilde B^+_{n-1} = -  \tilde B^+_{n+1}$.
Thus $k$ is determined by the coefficients $\tilde C^{\pm}_n$.
If $k$ is even and not a multiple of $m$ then
$\tilde B^+_{k} = -  \tilde B^+_{k}\ne 0$ and is minimal for 
this property.  
\end{proof}

\begin{example}\label{example3}
Let $f(x)=x^m$, $m$ odd, and let $g_1 (y) =  y^{km}$,
$g_2 (y) =  - y^{km}$, $k$ even.
The total zeta functions of $g_1$ and $g_2$ are equal but
the positive and the negative ones are different.
By Thom-Sebastiani formulae \eqref{tildeC} and example 
\ref{example2} all zeta functions of $f*g_1$ and $f*g_2$  
 coincide.  The functions $f*g_1$ and $f*g_2$ are not analytically 
equivalent but we shall show in the proof of theorem 
\ref{classification} below that they 
are blow-analytically equivalent.  
\end{example}

For the proof of theorem \ref{Thom-Seb2} we need the following lemmas.

\begin{lem}\label{tildeAn}
\begin{eqnarray*}
& \tilde{A}_n^+ & = (-1)^{nd_1} \chic (\{ \gamma \in \alphan;
f \circ \gamma = ct^n + \cdots , c \ge 0 \}) \\
& \tilde{A}_n^- & = (-1)^{nd_1} \chic (\{ \gamma \in \alphan;
f \circ \gamma = ct^n + \cdots , c \le 0 \}).
\end{eqnarray*}
\end{lem}

\begin{proof}
Denote by 
$$
\pi_{n,i} : \alphan \to \alphai
$$
the truncation map. It is a trivial fibration with fiber
isomorphic to
$\R^{(n-i)d_1}$ .  Then
$$
\alphan = \pi_{n,1} \inv (\ix_1) \cup
\pi_{n,2} \inv (\ix_2) \cup \cdots \cup
\pi_{n,n-1} \inv (\ix_{n-1}) \cup \ixnm \cup
\bixnp
$$
where by $\bixnp$ we denote $\{ \gamma \in \alphan;
f \circ \gamma = ct^n + \cdots , c \ge 0 \}$.  Then 
\begin{eqnarray*}
& \chic (\bixnp) &=  \chic (\alphan ) -
\sum_1^{n-1} (-1)^{(n-i)d_1} \chic (\ix_i) - 
\chic (\ixnm ) \hfill \\ 
&  &= (-1)^{nd_1} - \sum_1^{n-1} (-1)^{nd_1} a_i - 
(-1)^{nd_1} a^-_n = (-1)^{nd_1} (A_n +a_n^+).  
\end{eqnarray*}  
\end{proof}

\begin{lem}\label{trivial}
Let the map
$$
\varphi_i :\{ \gamma \in \alphan ; \ord f\circ \gamma
=i \} \to \R
$$
associate to $\gamma$, such that $(f\circ \gamma) (t) =
v_it^i + v_{i+1} t^{i+1} + \cdots$, the coefficient $v_n$.
Then, for $i<n$, $\varphi_i $ is a trivial fibration.  
\end{lem}

\begin{proof}
Define on $\{ \gamma \in \alphan ; \ord f\circ \gamma
=i \}$ an action of $\R$ by
$$
(\alpha, \gamma (t)) \to \gamma (t + \alpha t^{n-i+1}) .
$$
Then
$$
(f\circ \gamma) (t + \alpha t^{n-i+1}) =
v_it^i + \cdots + v_{n-1} t^{n-1} + (v_n +i v_i \alpha ) t^n + \cdots
$$
that gives
$$
\varphi_i (\gamma (t + \alpha t^{n-i+1})) = 
\varphi_i (\gamma ) + i v_i \alpha .
$$
Thus this action of $\R$ trivializes $\varphi_i$.  
\end{proof}

\smallskip
\begin{proof}[Proof of theorem \ref{Thom-Seb2}]
We show the formula for $\tilde C^+_n$.   By lemma \ref{tildeAn}
$$
\tilde C^+_n = (-1)^{nd} \chic
( \bixnp (f*g)), 
$$
where $\bixnp (f*g) = \{ (\gamma_1,\gamma_2) \in (\alphan (f) 
\times \alphan (g))  ; f(\gamma_1(t)) + g(\gamma_2(t))
= ct^n + \cdots, c \ge 0 \} $.  Then either
$\ord f(\gamma_1(t)) \ge n$ and $\ord g(\gamma_2(t))\ge n$
or $\ord f(\gamma_1(t)) =  \ord g(\gamma_2(t)) <n$.  
This gives the following decomposition
\begin{equation}\label{decomp}
\bixnp (f*g) = (Z\cap \bixnp (f*g)) \cup \bigcup_{i=1}^{n-1} 
(Z_i \cap \bixnp (f*g)) ,
\end{equation}
where
$$
Z= \{ (\gamma_1,\gamma_2)  ; \ord f(\gamma_1(t)) \ge n,  
\ord g(\gamma_2(t))\ge n\}
$$
and
$$
Z_i = \{ (\gamma_1,\gamma_2)  ; \ord f(\gamma_1(t))=  
\ord g(\gamma_2(t))=i\}.
$$

First we shall compute  $\chic (Z\cap \bixnp (f*g))$. 
Consider the map
$$
\Phi : Z \to \R^2_{(\varphi,\psi)}
$$
that associates to the pair of arcs 
$(\gamma_1,\gamma_2) $  the coefficients at 
$t^n$ of $f(\gamma_1(t))$ and of $g(\gamma_2(t))$.  
$\Phi$ is trivial over the following strata of
$\Phi(Z\cap \bixnp (f*g)) \subset \R^2$: 
$\{\varphi >0 , \psi>0 \}$, 
$\{\varphi + \psi>0 , \psi<0 \}$, 
$\{\varphi + \psi>0 , \varphi<0 \}$,
$\{\varphi >0 , \psi=0 \}$, 
$\{\psi>0 , \varphi=0 \}$,
$\{\varphi + \psi=0 , \psi<0 \}$, 
$\{\varphi + \psi=0 , \varphi<0 \}$,
and $\{\varphi = \psi=0 \}$. 
Note that $\Phi$ is trivial over
$\{\varphi + \psi \ge 0 , \psi<0 \}$ and
$\{\varphi + \psi \ge 0 , \varphi<0 \}$.
By this triviality  
$$
\chic (\Phi \inv (\{\varphi + \psi \ge 0,  
\psi<0 \})) =0
$$
since  $\chic (\{(\varphi,\psi)\in \R^2; 
\varphi + \psi \ge 0, \psi<0 \}) =0$.  Similarly 
$$
\chic (\Phi \inv (\{\varphi + \psi \ge 0, \varphi <0, \})) =0. 
$$
 Hence, since 
$Z\cap \bixnp (f*g) = \Phi \inv (\{\varphi + \psi \ge 0 \})$, 
\begin{eqnarray}\label{chicZ}
& \chic (Z\cap \bixnp (f*g))&
= \chic (\Phi \inv (\{\varphi + \psi \ge 0 \}))
=\chic (\Phi \inv (\{\varphi \ge 0,  \psi \ge 0 \})) \\
& & =(-1)^{nd_1} \tilde A_n^+ (-1)^{nd_2} \tilde B_n^+ 
=(-1)^{nd} \tilde A_n^+ \tilde B_n^+. \nonumber 
\end{eqnarray}

Now we compute 
$\chic (Z_i\cap \bixnp (f*g))$ for $0<i<n$ fixed.
Let $(\gamma_1,\gamma_2) \in Z_i\cap \bixnp (f*g)$.
Write
\begin{equation}\label{diei}
f(\gamma_1(t)) = v_it^i + \cdots + v_nt^n +\cdots , \quad 
 g(\gamma_2(t)) = w_it^i + \cdots + w_nt^n +\cdots . 
\end{equation}
Then $v_i\ne 0$, $w_i\ne 0$, and $v_i+w_i =0$.
Thus $v_i$ and $w_i$ are of opposite signs and hence 
$Z_i\cap \bixnp (f*g)$ is the disjoint union of two sets 
$$
Z_i^\pm = (\pi \inv _{n,i} (\ix _{i,\pm}(f)) \times
\pi \inv _{n,i} (\ix _{i,\mp}(g)))\cap \bixnp (f*g), 
$$
where $\pm$ equals the sign of $v_i$.  
Consider the following map 
$$
\Psi : \pi \inv _{n,i} (\ix _{i,+}(f)) \times
\pi \inv _{n,i} (\ix _{i,-}(g)) \to
\R_{(\varphi,\psi)}^{2}
$$
that associates to  $(\gamma_1,\gamma_2)$ the 
coefficients $(v_n,w_n)$ of \eqref{diei}.
By lemma \ref{trivial}, $\Psi$ is a trivial fibration and, 
since $\chic (\{(\varphi,\psi)\in \R^2; 
\varphi + \psi \ge 0 \}) =0$, 
$$
\chic (Z_i^+ )  = \chic (\Psi \inv (\{ 
\varphi + \psi \ge 0 \})) = 0. 
$$
Similarly we show that  $\chic (Z_i^-)=0$ and hence 
\begin{equation}\label{chicZi}
\chic (Z_i\cap \bixnp (f*g)) = 0.  
\end{equation}
The required formula for $\tilde C_n^+$ now follows from
\eqref{decomp}, \eqref{chicZ}, \eqref{chicZi}.  
\end{proof}

The formulae of theorem \ref{Thom-Seb} and the ones of theorem \ref{Thom-Seb2} 
are equivalent that one may check easily by a long but elementary computation. 
Alternatively, theorem \ref{Thom-Seb} can be proved by a topological argument 
similar to that of the proof of theorem \ref{Thom-Seb2}.  
We sketch just the main steps below.  
The details are left to the reader.  

\begin{proof}[Proof of theorem \ref{Thom-Seb}]
First note that the proof of lemma \ref{tildeAn} gives also  
$$
A_n = (-1)^{nd_1} \chic (\{ \gamma \in \alphan;
\ord f\circ \gamma >n\}) .
$$
Then $(\alphan (f) \times \alphan (g)) \cap \ixnp (f*g) =
(Z\cap \ixnp (f*g)) \cup \bigcup_{i=1}^{n-1} (Z_i \cap
\ixnp (f*g))$,  
with $Z$ and $Z_i$ as before.  By the triviality of 
$\Phi : Z \to \R^2_{(\varphi,\psi)}$ over the strata we get 
\begin{equation*}
(-1)^{nd} \chic (Z\cap \{\varphi + \psi >0\}) =
 a_n^+b_n^+ + a_n^+ b_n^- + a_n^- b_n^+  +
A_n b_n^+ +  a_n^+ B_n .  
\end{equation*}
Another argument based on lemma \ref{trivial} gives 
\begin{equation*}
(-1)^{nd} \chic (Z_i\cap \ixnp (f*g)) =
(-1)^{n-i} (a_i^+ b_i^- + a_i^- b_i^+) .  
\end{equation*}
Formula \eqref{c+} now follows from the additivity of Euler 
characteristic with compact supports.   
\end{proof}

%%%%%%%%%%%%%%%%%%%%%%%%%%%%%%%%%%%%%%%%%%%%%%%%%%%%%%%%%%%%%%%%%%%%

\bigskip
\section{Computations in two variables case}
\label{Toric}
\medskip

In this section we compute two dimensional examples using 
toric resolution.  
First we recall briefly the construction of
 toric resolution associated to a system of weights. 

Given a weight vector $(m,k)\in \N^2$, $m$ and $k$ coprime.
There is a canonical decomposition of the closed first
quadrant $\R_{\ge 0} \times \R_{\ge 0}$  in $\R^2$ into
a finite union of nonsingular
rational convex polyhedral cones that is
compatible with the weight vector.

\epsfxsize=6cm
\epsfysize=3.3cm
$$\epsfbox{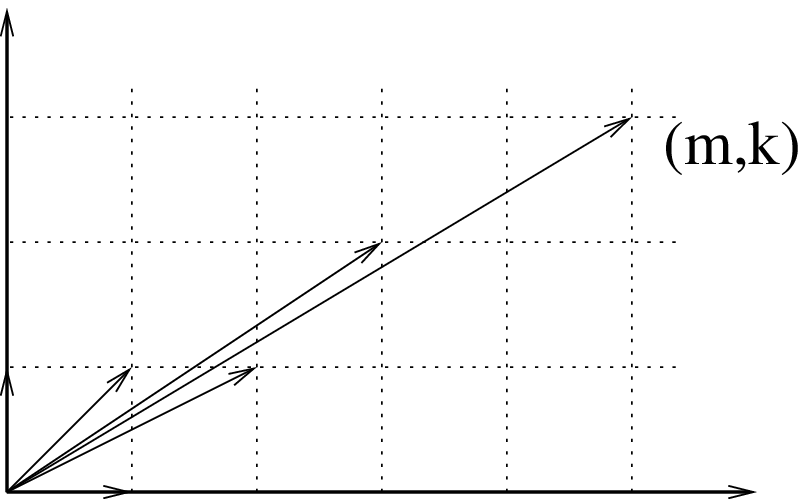}$$ 
%\vglue5cm

This decomposition
induces a toric modification $\sigma: M_\Delta \to \R^2$,
where $\Delta $ is the fan associated to this decomposition
and $M_\Delta$ is the associated toric variety. The exceptional
divisors of $\sigma$ are in one-to-one correspondence with the
one dimensional subcones (called rays or edges) of
$\Delta$ that are not the coordinate half-axis.
The integral vectors that generate these rays can be computed
out of $m,k$ by the following procedure.  Consider the
Hirzebruch-Jung continued fraction of $\frac m k $
$$
\frac m k = a_1 - \frac 1 {\displaystyle
{a_2 -
\frac 1 {\displaystyle { \cdots - \frac 1 {a_r} }}}}, 
$$
where $a_i \ge 2$ for $i>1$ and $a_1\ge 1$.  The
coefficients $a_i$  define the 
vectors $(m_i,k_i)\in \R^2$, $i=1, \cdots, r+1,$ such that 
\begin{equation*}
\begin{array}{cccc}
m_1 = 1,  & m_2 = a_1, &  m_{i+1} = a_i m_i - m_{i-1}  &
\text {for } 2\le i\le r   ;   \\
k_1 =0,   & k_2 = 1,   & k_{i+1} = a_i k_i - k_{i-1}   &
 \text {for } 2\le i\le r  , 
\end{array}
\end{equation*}
and then $m_{r+1}=m$, $k_{r+1}=k$.  
Similarly the coefficients $b_1, \cdots, b_s$ of the
Hirzebruch-Jung continued fraction of $\frac k m$ 
define the 
vectors $(m'_i,k'_i)\in \R^2$, $i=1, \cdots, s+1,$ 
such that 
\begin{equation*}
\begin{array}{cccc}
m'_1 = 0,  & m'_2 = 1, &  m'_{i+1} = b_i m'_i - m'_{i-1}  &
\text {for } 2\le i\le s   ;   \\
k'_1 =1,   & k'_2 = b_1,   & k'_{i+1} = b_i k'_i - k'_{i-1}   &
 \text {for } 2\le i\le s   
\end{array}
\end{equation*}

Then the vectors 
\begin{equation}\label{vectors}
(1,0) = (m_1,k_1), \ldots, (m_r,k_r), (m,k), (m'_s,k'_s),
\ldots ,  (m'_1,k'_1)=(0,1) 
\end{equation}
are the primitive vectors of the rays of $\Delta$.   Choose a pair 
 of subsequent vectors $\mathbf v = (a,b)$, $\mathbf w=  (c,d)$ 
of \eqref{vectors}.  They
generate a two dimensional cone $\tau$ of $\Delta$ and give rise 
to an affine chart of $\sigma$, 
$\sigma_\tau: M_\tau \simeq \R^2 \to \R^2$ given by 
$$
\sigma_\tau (X,Y) = (X^aY^c,X^bY^d) . 
$$ 
The divisor corresponding to $\mathbf v$, resp. $\mathbf w$, is given 
in $M_\tau$ by $X=0$, resp. $Y=0$.  
The jacobian 
$$
\jac \sigma_\tau = \left | 
\begin{array}{cc}
a  & c   \\
b   & d   
\end{array}
\right | X^{a+b-1}Y^{c+d-1} = X^{a+b-1}Y^{c+d-1}  
$$
and hence it is normal crossings.  Denote  by $E_\mathbf v$
the divisor corresponding 
to the vector $\mathbf v$.  Then the multiplicity of 
$\jac \sigma_\tau$ along $E_\mathbf v$ equals
$$
\mult_{E_\mathbf v} \jac \sigma_\tau = a + b -1 . 
$$  
Let 
$$
f(x,y) = \sum_{(i,j)\in \N^2 \setminus \{(0,0)\}} 
a_{i,j} x^i y^j 
$$
and denote $\supp (f) = \{(i,j); a_{i,j}\ne 0\}$.   
Then 
$$
\mult_{E_{\mathbf v}} f\circ \sigma =
\min_{(i,j)\in \supp (f)} \{ai+bj\} .
$$ 
 
%\medskip
\begin{example}\label{example4}
We compute the toric resolution and the zeta functions of 
$f(x,y) = x^3 + xy^5$.  $f$ is nondegenerate weighted homogeneous
 with weights $(5,2)$.  The toric modification associated to 
this system of weights is given by the 
vectors $(1,0)$, $\mathbf v_1 = (3,1)$, $\mathbf v_2 =(5,2)$, 
$\mathbf v_3=(2,1)$, $\mathbf v_4= (1,1)$, $(0,1)$.  

\epsfxsize=5.5cm
\epsfysize=2.5cm
$$\epsfbox{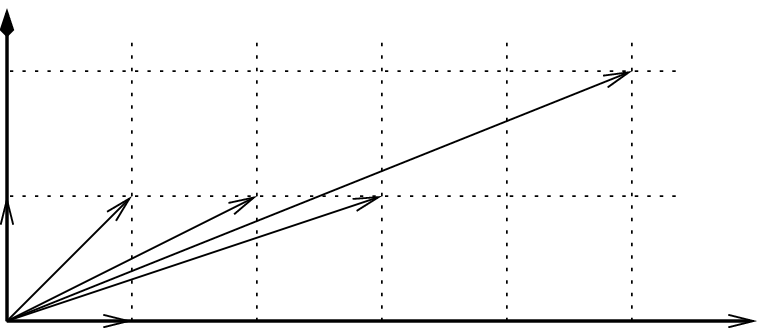}$$
%\vglue4.5cm

Denote by $E_i$ the component of the exceptional divisor corresponding
to $\mathbf v_i$.  Let $N_i = \mult_{E_i} f\circ \sigma$, 
$\nu_i = \mult_{E_i} \jac \sigma +1$.   Then, by above, $N_1 = 8$, 
$\nu_1 = 4$, $N_2 = 15$, $\nu_2 = 7$, $N_3 = 6$, 
$\nu_3 = 3$, $N_4 = 3$, $\nu_4 = 2$.  The strict transform 
of the zero set of $f$ has two components;
the strict transform $S_1$ of 
$x=0$ and the strict transform $S_2$ of $x^2 + y^5=0$.
The first one intersects $E_1$ and the second one $E_2$
as indicated on the resolution tree below.   

\epsfxsize=9.5cm
\epsfysize=4cm
$$\epsfbox{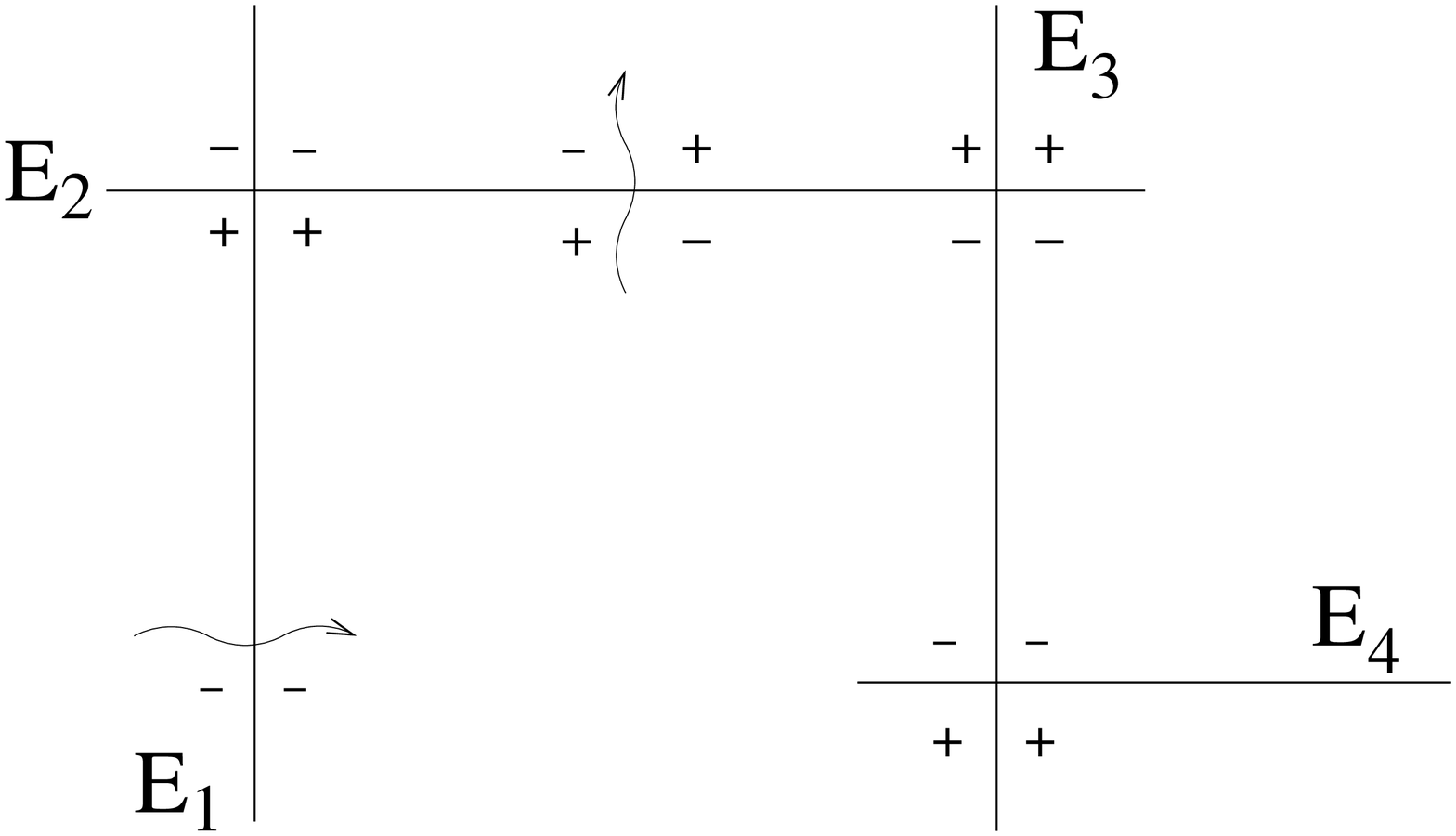}$$
%\vglue4cm

Thus 
\begin{eqnarray}\label{xy5} 
& & Z(T)  =  4 \frac { T^{8}} {1- T^{8}} - 6 \frac { T^{15}} {1+ T^{15}}
- 4 \frac { T^{6}} {1+ T^{6}} + 2 \frac { T^{3}} {1-  T^{3}} 
- 4  \frac { T^{8}} {1- T^{8}} \frac { T^{15}} {1+ T^{15}} \\
& & \quad 
 +4 \frac { T^{15}} {1+ T^{15}}  \frac { T^{6}} {1+ T^{6}} 
- 4 \frac { T^{6}} {1+ T^{6}} \frac { T^{3}} {1-  T^{3}} 
-4 \frac { T^{8}} {1- T^{8}} \frac T {1+T} +4 \frac { T^{15}} 
 {1+ T^{15}}\frac T {1+T} .  \nonumber
\end{eqnarray}
and $Z_+(T) = Z_- (T)= \half Z(T)$.  
\end{example}

%%%%%%%%%%%%%%%%%%%%%%%%%%%%%%%%%%%%%%%%%%%%%%%%%%%%%%%%%%%%%%%%%

\bigskip
\section{Zeta functions are blow-analytic invariants}
\label{Zetainvariant}
\smallskip

Blow-analytic equivalence is a notion introduced
by T.-C. Kuo as a natural equivalence relation 
for real analytic function germs.
He established several fundamental results on blow-analyticity.
For a general review on the blow-analytic theory (until 1997),
see \cite {fukuikoikekuo}.
The notion of blow-analytic equivalence 
is defined as follows:

\vspace{3mm}

We say that analytic function germs 
$f, g : (\R^d,0) \to (\R,0)$ are {\em blow-analytically equivalent} 
if there are real modifications
$\mu : (M,\mu^{-1}(0)) \to (\R^d,0)$,
$\mu^{\prime} : (M^{\prime},\mu^{\prime -1}(0)) \to (\R^d,0)$
and an analytic isomorphism $\Phi : (M,\mu^{-1}(0)) \to
(M^{\prime},\mu^{\prime -1}(0))$
which induces a homeomorphism 
$\phi : (\R^d,0) \to (\R^d,0)$
such that $f = g \circ \phi$.

\vspace{3mm}

By a real modification, we mean the following.
Let $\mu : M \to N$ be a proper surjective analytic map 
of real manifolds.
It has a unique extension to a holomorphic map
$\mu^* : U(M) \to U(N)$ where $U(M), \ U(N)$ are
respectively open neighborhoods of $M, \ N$
in their complexifications $M^*, \ N^*$.
We say that $\mu$ is a {\em real modification}
if $\mu^*$ is an isomorphism except on some 
thin subset of $U(M)$. 

Let $\mu : (M,\mu^{-1}(0)) \to (\R^d,0)$ be a real modification.
Take any analytic arc at $0 \in \R^d$,
$\lambda : (-\epsilon , \epsilon ) \to \R^d$,
$\lambda (0) = 0$.
Then $\lambda$ has an analytic lifting.
Namely, there is an analytic arc 
$\lambda^{\prime} : (-\epsilon , \epsilon ) \to M$,
$\lambda^{\prime} (0) = P \in \mu^{-1}(0)$
such that $\lambda^{\prime} \circ \mu = \lambda$.
Remark that if $\lambda$ is not contained 
in the critical value set of $\mu$ (a thin subset of $\R^d$)
as set-germs at $0 \in \R^d$, then the lifting is unique.

In this paper, we assume also the following condition
for the real modifications $\mu$ and $\mu^{\prime}$
in the definition of blow-analytic equivalence: 
the critical value sets of $\mu$ and $\mu^{\prime}$
are contained in the zero-sets of $f$ and $g$
respectively as set-germs at $0 \in \R^d$.
The assumption is reasonable.
In fact, for any analytic function germ $f : (\R^d,0) \to (\R,0)$, 
there is a real modification $\mu : (M,\mu^{-1}(0)) \to (\R^d,0)$ 
with this property such that $f \circ \mu$ is a normal crossing 
(\cite{hironaka}, \cite{bierstonemilman}).
Any triviality theorem (\cite{kuo1}, \cite{fukuiyoshinaga},
\cite{fukuipaunescu}, \cite{abderrahmane} and so on) 
and a locally finite classification theorem
(\cite{kuo2}) have been established on blow-analytic
equivalence with the property.
A blow-analytic invariant (e.g. \cite{fukui}) in the original sense 
is, of course, a blow-analytic invariant in our sense.

Suppose that real analytic function germs
$f, g : (\R^d,0) \to (\R,0)$ are blow-analytically equivalent
in the sense of this paper.
Then we can say that
the uniqueness of the arc lifting property holds
for $\mu$ (resp. $\mu^{\prime}$) 
if the arc is not contained in a subset of
the zero-set $f^{-1}(0)$ (resp. $g^{-1}(0)$).
   
\smallskip
In this section we show that if two analytic function germs 
$f, g :(\R^d,0) \to (\R, 0)$ are blow-analytically  equivalent 
then their zeta functions 
%co{\"\i}ncide 
coincide that is $\zetafa = Z_g, 
\zetafp = Z_{g,+}, \zetafm = Z_{g,-}$.  This will follow 
from Denef \& Loeser's formulae \eqref{DL}, \eqref{DLpm} 
that we show first.  
The proof will be an adaptation to the real analytic geometry,
the ideas of \cite{kontsevich}, \cite{denefloeser2}.  
The main difficulty is that 
%the set 
we have to use the sets that are 
not necessarily semi-algebraic 
but only subanalytic and not relatively compact, 
so we have to show that they have a well defined 
Euler characteristic with compact supports that is additive.  

Let $\sigma:(M,E_0) \to (\R^d,0)$, $E_0=\sigma \inv (0)$, be a real 
modification.  
%M analytic manifold, \mu analytic proper, \mu satisfies the arc-lifting 
Consider the space of analytic arcs   
$$
\alphaanal (M,E_0) 
:= \{ \gamma: (\R,0)\to (M,E_0); \gamma 
\text { analytic} \} .  
$$
The set of truncated arcs can be described  
as follows 
 $$
\alphan (M,E_0):= \alphaanal (M,E_0) /\sim ,
$$
where $\gamma_1(t) \sim \gamma _2 (t)$ if $\gamma_1(0) = \gamma_2 (0)$
and $\gamma_1 (t) - \gamma_2(t) = O(t^{n+1})$ in a (or any) 
local system of coordinates at $\gamma_1(0) = \gamma_2 (0)$.  
$\alphan (M,E_0)$ is an analytic variety, a subvariety of a similarly
defined set $\alphan (M)$ that is an analytic manifold.  The projection 
$\alphan (M,E_0)\to E_0$ is a locally trivial fibration with fiber
$\R^{nd}$. 
Indeed, in a local system of coordinates on an open neighborhood
$U$ of $p_0\in E_0$ we may write simply 
\begin{equation}\label{y}
\alphan (U,U\cap E_0) = \{ \gamma; \gamma(t) = p +
\mathbf y_1 t + \cdots + \mathbf y_n t^n,  p\in U\cap E_0, 
\mathbf y_i\in \R^d \}.  
\end{equation}
Denote  $\mathbf y := (\mathbf y_1,\ldots ,\mathbf y_n)$.
Using the coordinates $(p,\mathbf y )$ we identify 
$\alphan (U,U\cap  E_0) \simeq (U\cap E_0) \times \R^{nd}$.  
Following {\L}ojasiewicz \cite{lojasiewicz} we call a semi-analytic
subset of $\alphan (U,U\cap  E_0)$ {\it relatively semi-algebraic 
with respect to} $\mathbf y$ if it is defined 
by a finite number of equations and inequalities in functions that are 
analytic in $p$ and polynomial in $\mathbf y$.  
Let us compare two such trivializations.  
This amounts to consider the following situation.  
Let $U,U'$ be two open subsets of $\R^d$ and let 
$h:U\to U'$ be an analytic isomorphism. 
Let $\gamma(t) = p + \mathbf y_1 t + \cdots + \mathbf y_n t^n$ 
be as above. Then 
$$
h(\gamma(t)) = h(p) + \va _1 (p,\mathbf y) t +  \cdots +  
\va _n (p,\mathbf y) t^n  + O(t^{n+1}).
$$
The coefficients $ \va _i (p,\mathbf y)$ are analytic in 
$p$ and polynomial in $\mathbf y$.  
Thus two such local trivializations of 
$\alphan (M,E_0) \to E_0$ differ by an analytic 
isomorphism that is polynomial on the fibers.  
%\begin{defn}
A semi-analytic subset $A$ of  $\alphan (M,E_0)$ will be called 
{\it relatively semi-algebraic} 
if for each $p\in E_0$ there is an open neighborhood 
$U$ of $p$ in $M$ such that  $A \cap \alphan (U,U\cap  E_0)$ is
relatively semi-algebraic.  
%\end{defn} 

Let $X$ be an analytic manifold.  If $A\subset X$ is subanalytic and
relatively compact then its (co)homology groups are finitely
generated.  The Euler characteristic (standard or with compact
supports) of such sets is well defined and the Euler characteristic
with compact supports is additive.  
This is not, in general, true if  $A$ is no longer relatively compact.  
This observation justifies the following definition. 

\begin{defn}
Let $X$ be an analytic manifold and $A\subset X$.  We say that  $A$ is 
{\it globally subanalytic} if there is an analytic manifold $\tilde X$
and an analytic embedding $i:X\to \tilde X$  such that 
$i(A)$ is relatively compact and subanalytic in $\tilde X$.  
\end{defn}

A trivial example of globally subanalytic sets are semi-algebraic 
subsets of $\R^N$.  The example we really have in mind are the
semi-analytic and relatively semi-algebraic subsets of
$\alphan (M,E_0) \to E_0$.   Indeed, we may suppose that $M$ is a
submanifold of $\R^N$ and hence $\alphan (M,E_0) \overset i
\hookrightarrow \R^N \times \R^{nN}$.  Choose any algebraic
compactification of $\R^{nN}$, the one point  compactification
$S^{nN}$ for instance.  If $A\subset \alphan (M,E_0)$ is
semi-analytic and relatively semi-algebraic then $i(A)$ is relatively 
compact and subanalytic (even semi-analytic) in $\R^N \times S^{nN}$.

Let $\pi_n:\alphaanal (M,E_0) \to \alphan (M,E_0)$ and $\pi_n:
\alphaanal (\R^d,0) \to \alphan (\R^d,0)$ denote the standard 
projections.  
The real modification $\sigma : (M,E_0) \to (\R^d,0)$ induces a map 
$$
\sigma_*: \alphaanal (M,E_0) \to 
\alphaanal (\R^d,0), 
$$ 
defined by composition $\sigma_*(\gamma)(t) = \sigma(\gamma(t))$ 
that gives an analytic map on truncations 
$$
\sigma_{*n}: \alphan (M,E_0) \to 
\alphan (\R^d,0).    
$$ 
Clearly $\pi_n \circ \sigma_* = \sigma_{*n} \circ \pi_n$.  

Let $\gamma \in \alphaanal (M,E_0)$.  The jacobian determinant 
$\jac \sigma$ of $\sigma$ may be defined using any local coordinate
system on $M$.  Its order in $t$ along $\gamma (t)$,  
$\ord \jac \sigma(\gamma (t))$,  
is independent of this choice.   
Given a positive integer $e$.  Define   
$\Delta_e = \{ \gamma \in \alphaanal (M,E_0); 
\ord \jac \sigma(\gamma (t)) = e \} $ and $\Delta_{e,n} = 
\pi_n(\Delta_e)$.

\medskip
\begin{lem}\label{keylemma}
Let $e\ge 1$ and $n\ge 2e$.  
\begin{enumerate}
\item [(a)] 
Let $\gamma_1,\gamma_2 \in \alphaanal (M,E_0)$. 
If $\gamma_1 \in \Delta_e$ and 
$\sigma (\gamma_1) \equiv \sigma (\gamma_2) \mod t^{n+1} $ then 
$\gamma_1 \equiv \gamma_2 \mod t^{n+1-e}$ and $\gamma_2 \in \Delta_e$. 
\item [(b)]
$\sigma_{*n} (\Delta_{e,n})$ is a globally subanalytic subset of 
$\alphan (\R^d,0)$. There exists a subanalytic stratification of 
$\sigma_{*n} (\Delta_{e,n})$ such that over each stratum
$\sigma_{*n}$ is a trivial fibration with fiber $\R^e$.  
\end{enumerate}
\end{lem}

\begin{proof}
Let $p\in E_0$.  Choosing a local coordinate system at 
$p$ we may suppose that $p=0\in \R^d$.  Let $\gamma (t)\in \Delta_e$, 
$\gamma(0) = 0$.  Denote by $\Jac_\sigma (x)$ the jacobian matrix of 
$\sigma$ at $x$.  
Then 
$$
\M (t) := t^e (\Jac_\sigma (\gamma (t)))\inv 
$$
is a matrix with entries analytic functions in $t$.  
By Taylor formula
\begin{equation}\label {taylor}
\sigma (\gamma (t) + t^{n+1-e} u) = \sigma (\gamma(t)) + 
t^{n+1-e} \Jac_{\sigma}(\gamma) u + R(\gamma (t),u) ,  
\end{equation}
where $R(\gamma (t),u)$ is analytic in $t$ and $u\in \R^d$.
Moreover, $R(\gamma (t),u)$ is divisible by $t^{2(n+1-e)}$ 
and hence by $t^{n+2}$. Let 
$$
R(\gamma (t),u) = t^{n+2} \tilde R(\gamma (t),u), \quad 
\tilde R(\gamma (t),u) \text{ analytic}. 
$$
We solve the following equation with respect 
to $u\in \R^d$
\begin{equation}\label{equation}
\sigma (\gamma (t) + t^{n+1-e} u) = \sigma (\gamma (t)) + t^{n+1} v.
\end{equation}
By \eqref{taylor}, \eqref{equation} is equivalent to 
\begin{equation*}%\label {taylor}
t^{n+1} v = 
t^{n+1-e} \Jac_\sigma (\gamma) u + R(\gamma (t),u) ,  
\end{equation*}
and hence to 
\begin{equation*}%\label {taylor}
u = \M (t) v - t \M (t) \tilde R(\gamma (t),u) .    
\end{equation*} 
By the Implicit Function Theorem,
for any $v_0 \in \R^d$, 
this equation has a unique analytic 
solution $u=u(t,v)$ defined in a neighborhood of $(v_0,0)$.  
In particular, if $v(t)$ is an analytic arc 
then \eqref{equation} admits a solution 
being an analytic arc $u(t)=u(t,v(t))$.  
This solution is unique since $\sigma$ is 
a modification (and $\jac \sigma (\gamma (t) + t^{n+1-e} u(t) )$ 
is not identically equal to $0$, see \eqref{taylor2} below).  

Now we show (a).  
Let $\sigma (\gamma_1) \equiv \sigma (\gamma_2) \mod t^{n+1} $ and 
consider a local coordinate system at $p = \gamma_1(0) = \gamma_2(0)$. 
By the above $\gamma_2(t)$ as the solution of \eqref{equation} with 
$\gamma = \gamma_1$ and $v(t) = t^{-(n+1)} (\sigma(\gamma_2)(t) - 
\sigma(\gamma_1)(t))$ 
is of the form 
$$
\gamma_2(t) = \gamma_1(t) + t^{n+1-e} u(t) .
$$
This shows the first claim of (a).  
By Taylor formula 
\begin{eqnarray}\label {taylor2}
\ \ \ \jac \sigma (\gamma_1 (t) + t^{n+1-e} u(t) )  
& = & \jac \sigma (\gamma_1(t)) + 
t^{n+1-e} \Jac_{\jac \sigma} (\gamma_1) u (t) + O(t^{2(n+1-e)}) \\
 & \equiv & \jac \sigma (\gamma_1(t)) \mod t^{e+1} , \hfill  \nonumber
\end{eqnarray}
since $n+1-e \ge e+1$.  This completes the proof of (a).   

We show (b).  
By (a) the set $\Delta_{e,n}$ is the union of fibers of $\sigma_{*n}$.  
To compute these fibers we fix $\gamma (t)\in \Delta_e$.  
We keep the notation of the first part of proof of lemma.      
By \eqref{equation}, the fiber of 
$\sigma_{*n}$ over $\pi_n (\sigma_* (\gamma))$ equals   
\begin{eqnarray*}
 \sigma_{*n}\inv (\pi_n (\sigma_* (\gamma)))  & = &\{ \gamma (t) +
t^{n+1-e} u \mod t^{n+1} ; u= \mathbf u_0 + \mathbf u_1 t + \cdots +  
\mathbf u_{e-1} t^{e-1},    \\
 &  & \quad \Jac _\sigma (\gamma (t)) u(t) \equiv 0 \mod t^e \} \hfill
\end{eqnarray*} 
and hence is isomorphic to a linear subspace of 
$\{ u= \mathbf u_0 + \mathbf u_1 t + \cdots +  
\mathbf u_{e-1} t^{e-1}; \mathbf u_i\in \R^d\} \simeq \R^{de}$.  
There are invertible matrices $A$ and $B$  with
entries in $\R\{t\}$ such that $A \Jac _\sigma (\gamma (t)) B$ is
equivalent over $\R\{t\}$ to a diagonal 
matrix with diagonal elements $t^{e_1}, \ldots, t^{e_d}$. 
(For this it suffices to apply to 
$\Jac _\sigma (\gamma (t))$ Gauss' elimination method.)
Necessarily $e = e_1 + \cdots + e_d$ and hence the fiber is isomorphic
to $\R^e$.  

The map $\sigma_{*n}: \alphan (M,E_0) \to \alphan (\R^d,0)$ 
is analytic but not proper.  
Therefore, even if $\Delta_{e,n}$ is a semi-analytic set, 
it is not immediate that its image $\sigma_{*n} (\Delta_{e,n})$ 
is subanalytic.  
This follows from the relative semi-algebraicity of 
$\sigma_{*n}$ and $\Delta_{e,n}$.  
By this we mean the following.  
Let $\Gamma \subset \alphan (M,E_0)\times \alphan (\R^d,0)$ be the
graph of $\sigma_{*n}$. 
Using a local system of coordinates at $p_0 \in E_0$ we
identify an open neighborhood of $p_0$ in $M$ with 
an open neighborhood $U$ of the origin in $\R^d$ 
so that $p_0$ corresponds to the origin.  
Then ${\sigma_{*n}}$ restricted to $\alphan (U,U\cap E_0)$ 
can be computed as follows.  Write $y(t) \in \alphan (U,U\cap E_0)$
as in \eqref{y} and $x(t) \in \alphan (\R^d,0)$ as 
$$
x(t) = \mathbf x_1 t + \cdots + \mathbf x_n t^n,
\quad \mathbf x_i\in \R^d
$$
Denote $\mathbf x := (\mathbf x_1,\ldots ,\mathbf x_n)$.  
Then each coefficient of $x(t) = \sigma_{*n} (y(t))$, 
$\mathbf x_j (p, \mathbf y)$ is analytic in $p$ 
and polynomial in $\mathbf y$. 
That is in these coordinates $\Gamma$ is given by an 
analytic equation $\mathbf x - \sigma_{*n} (p, \mathbf y)=0$ 
that is polynomial 
in $(\mathbf y,\mathbf x)$.
We say for short that  $\Gamma$ is relatively semi-algebraic 
with respect to the projection onto $E_0$.  
Let $\Gamma_\Delta \subset \Gamma$ 
be the graph of $\sigma_{*n}$ restricted to  $\Delta_{e,n}$.
A similar argument  shows that $\Gamma_\Delta $ is relatively
semi-algebraic with respect to the projection onto $E_0$.  
Therefore, by {\L}ojasiewicz's version of Tarski-Seidenberg Theorem 
\cite{lojasiewicz}, the projection $pr (\Gamma_\Delta )$ of 
$\Gamma_\Delta$ into $E_0 \times \alphan (\R^d,0)$ is semi-analytic
and relatively semi-algebraic.  Finally, since $E_0$ is compact,
the projection of $pr(\Gamma_\Delta )$ in $\alphan (\R^d,0)$, that
equals $\sigma_{*n} (\Delta_{e,n})$, is subanalytic.  Moreover, it is 
easy to see that it is globally subanalytic.

Our original identification of $\alphan (\R^d,0)$ with the space of
truncated arcs $x(t) = \va _1 t + \va _2 t^2 + \cdots +  \va _n t^n$
gives an inclusion $\alphan (\R^d,0) \hookrightarrow \alphaanal
(\R^d,0)$ that is a section of \ $\pi_n$.  This allows us to define
a section $s$ of $\sigma _{n*}$ by 
\begin{equation}\label{section}
s: \alphan (\R^d,0) \hookrightarrow \alphaanal (\R^d,0) \overset 
{\sigma_*\inv} \to  
\alphaanal (M,E_0) \overset  {\pi_n} \to \alphan (M,E_0)
\end{equation}
that is defined on those curves that are not entirely contained 
in the crititcal locus of $\sigma$, 
in particular on $\sigma_{*n}(\Delta_{e,n})$.  
Let $s_\Delta$ be the restriction of $s$ onto 
$\sigma_{*n} (\Delta_{e,n})$ and let $\Gamma_{s,\Delta}$ be the graph 
of $s_\Delta$.  We shall show that $\Gamma_{s,\Delta}$ is 
globally subanalytic 
in $\alphan (\R^d,0) \times \alphan (M,E_0)$.  
Considering $ \sigma_{*n} (\Delta_{e,n}) $ as 
a subset of $\mathcal L_{n+e} (\R^d,0)$ by sequence of inclusions 
$\sigma_{*n} (\Delta_{e,n}) \subset 
\alphan (\R^d,0) \subset \mathcal L_{n+e} (\R^d,0)$ define 
$$
\tilde \Gamma_{s,\Delta} = \{(x(t), y(t)) \in 
\sigma_{*n} (\Delta_{e,n}) \times \mathcal L_{n+e}(M,E_0); 
x(t) = \sigma_{*(n+e)} (y(t))\}.  
$$
Then the graph $\Gamma_{s,\Delta}$  
of $s_\Delta$ is the projection of $\tilde \Gamma_{s,\Delta}$ to 
$\alphan (\R^d,0) \times \alphan (M,E_0)$.  
Indeed, it is clear that $\Gamma_{s,\Delta}$
is contained in this projection.
On the other hand, if $x(t) = \sigma_{*(n+e)} (y(t))$ then 
$\sigma_{*(n+e)} (y(t)) = \sigma_{*(n+e)} (s(x(t)))$ and  
by (a), $y(t) \equiv \hat s(x(t)) \mod t^{n+1}$.  
Note that $\tilde \Gamma_{s,\Delta}$ is a semi-analytic set
relatively semi-algebraic with respect to the projection to $E_0$,
and hence so is $\Gamma_{s,\Delta}$. 
Note also that $s_{\Delta}$ need not to be continuous and usually
it is not, see example \ref{illustrative} below.

Now we are ready to finish the proof of (b).  Fix a subanalytic 
stratification of $\sigma_{*n} (\Delta_{e,n})$ 
so that $s_{\Delta}$ is analytic on each stratum.  
Such a stratification exists since the graph of 
$s_{\Delta}$ is globally subanalytic. 
Fix a stratum $S$.  
Subdividing $S$ if necessary, we may suppose that $s(S)$ is contained
in an open subset of $\alphan (M,E_0)$ corresponding to a local chart
on $M$.  Thus we may use local coordinates on $M$.    
\begin{eqnarray*}
 \sigma_{*n}\inv (S) = \{ s(x)(t) + t^{n+1-e} u \mod t^{n+1} ;
x(t) \in S, \Jac _\sigma (s(x)(t)) u(t) \equiv 0 \mod t^e \}
\end{eqnarray*} 
Since the kernel of $\Jac _\sigma (s(x)(t))   \mod t^e$ is isomorphic 
to $\R^{nd}$ for each $x\in S$, $\sigma_{*n}$ is a locally trivial 
analytic fibration over $S$.  Thus subdividing again $S$, if necessary, 
we may ensure that the fibration becomes trivial over each stratum.  
This ends the proof.    
\end{proof}

Let $A\subset \alphaanal (M,E_0)$ (or  $A\subset \alphaanal (\R^d,0)$).  
We say that $A$ 
is {\it subanalytic} if $A= \pi_n \inv (C)$ where $C$ is a 
globally subanalytic 
subset of $\alphan (M,E_0)$ (resp. of $\alphan (\R^d,0)$).
We say that $A$ is  {\it $n$-stable} if $A$ is subanalytic and 
$A= \pi_n \inv (\pi_n(A))$. For instance by  Lemma \ref{keylemma}, 
$\sigma_* (\Delta_e)$ is  $2e$ stable. It follows from  
\eqref{taylor2} that   $\Delta_e$ is  always $e$ stable.

\begin{example}\label{illustrative}
$\sigma(X,Y)=(X^2Y,XY)$, $e=2$.

Let $\sigma:\R^2\to \R^2$ be given by $(x,y)=\sigma(x,y)= (X^2Y,XY)$. 
Consider the curve
\begin{eqnarray*}
 & & \hfil \gamma(t) = (X(t),Y(t)) = (X_0 + X_1t + X_2 t^2 + \cdots ,
Y_0 + Y_1t + Y_2 t^2 + \cdots) \\
& & \sigma(\gamma(t)) = (x(t),y(t)) = (x_0 + x_1t + x_2 t^2 + \cdots ,
y_0 + y_1t + y_2 t^2 + \cdots) . 
\end{eqnarray*}
The jacobian determinant $\jac \sigma = X^2Y$ and
\begin{eqnarray*}
& & \Delta_2  =  \{ \gamma(t) ; X_0Y_0=X_0Y_1 =0, X_0^2Y_2 + X_1^2 Y_0 
\ne 0 \} \\
& & \sigma_*(\Delta_2) =  \{\sigma(\gamma(t)); x_0=x_1=y_0=0, x_2\ne
0 \}.
\end{eqnarray*}
The conditions on $\Delta_2$ do not involve $X_i,Y_i$, $i>2$, so
$\Delta_2$ is 2-stable (as always).  In this example also
$\sigma_*(\Delta_2)$ is 2-stable.  The truncations $\Delta_{2,2}$
and $\sigma_{*2}(\Delta_{2,2})$, that are subsets of
$\alphaanal_2(\R^2) \simeq \R^6$, are given by the same
conditions.  Note that $\sigma_{*2}(\Delta_{2,2})$ is irreducible
but $\Delta_{2,2}$ has two irreducible components.  They are
2-truncations of
\begin{eqnarray*}
\Delta'_2 := \Delta_2 \cap \{X_0=0\}, \qquad
\Delta''_2 := \Delta_2 \cap \{Y_0=0\}, 
\end{eqnarray*}
and their images are respectively 
\begin{eqnarray*}
\sigma_*(\Delta'_2) = \sigma_*(\Delta_2) \cap \{y_1\ne 0\},  \qquad
\sigma_*(\Delta''_2) = \sigma_*(\Delta_2) \cap \{y_1=0\}.
\end{eqnarray*}
Thus by modification $\sigma$ the geometry of the set of curves in
$\sigma_*(\Delta_2)$ changes dramatically and $\Delta_2$ contains
the curves of two different kinds: the ones hitting $\{X=0, Y\ne
0\}$ transversally and the ones touching $\{Y=0, X\ne 0\}$ with
intersection number $2$.

Both restrictions of $\sigma_{*2}$: $\Delta'_2 \to
\sigma_*(\Delta'_2)$ and $\Delta''_2 \to
\sigma_*(\Delta''_2)$, are trivial fibrations with fiber $\R^2$.
For instance, the first one is given by
$$
x_0=x_1=y_0=0, \, x_2=X_1^2Y_0, \, y_1=X_1Y_0, \, y_2= X_1Y_1+ X_2Y_0.
$$
The section $s$ of $\sigma_{*2}$ is defined in \eqref{section}.  We 
compute the restriction of $s$ to $\sigma_{*2}(\Delta'_{2,2})$.  Fix
a curve in $\sigma_{*2}(\Delta'_{2,2})$
$$
(x(t),y(t)) = (x_2t^2, y_1t+y_2t^2)
$$
that we consider as a curve in $\sigma_{*}(\Delta'_{2})$.  It
lifts to
\begin{eqnarray*}
& & X(t)  = x(t) (y(t))\inv = x_2y_1\inv t (1- (y_2/y_1)t +\cdots ) \\
& & Y(t)  = (x(t))\inv (y(t))^2 = x_2\inv (y_1^2 + 2y_1y_2 t +
y_2^2 t^2).  
\end{eqnarray*}
That is $s$ on $\sigma_{*2}(\Delta'_{2,2})$ is given by 
$$
s(0,0,x_2,0,y_1,y_2)=(0,x_2/y_1, x_2y_1^{-2}y_2, x_2\inv y_1^2,
2x_2 \inv y_1y_2, x_2\inv y_2^2). 
$$
Recall that $x_2\ne 0$ everywhere on $\sigma_{*2}(\Delta_{2,2})$ but
$y_1$ vanishes  on $\sigma_{*2}(\Delta''_{2,2})$.  Thus $s$ cannot
be extended continuously from  $\sigma_{*2}(\Delta'_{2,2})$ to
$\sigma_{*2}(\Delta_{2,2})$.
A similar computation shows that $s$ on $\sigma_{*2}(\Delta''_{2,2})$
is given by $s(0,0,x_2,0,0,y_2)=(x_2/y_2,0,0,0,0, y_2^2/x_2)$.  
\end{example}

\medskip
By definition each subanalytic $A\subset  \alphaanal (M,E_0)$ is 
$n$ stable for $n$ sufficiently large.    
Following   
\cite{kontsevich}, \cite{denefloeser2}, \cite{denefloeser0}, we may 
associate to each $n$-stable $A$ its motivic measure that will be in
our case simply 
$$
\chic (A) := (-1)^{-(n+1)d} \chic (\pi_n (A)), .
$$
This expression is independent of $n$ (if $A$ is $n$-stable). 
We say that  $\varphi:A\to \Z$ is {\it constructible} if the 
image of $\varphi$ is finite 
and $\varphi \inv (m)$ is subanalytic for each $m\in \Z$.  Then we 
define 
$$
\int_A \varphi \, d\chic := \sum_{m\in \Z} m \chic (\varphi\inv (m)).
$$ 

The following corollary of Lemma \ref{keylemma} is a real analytic
version of Kontsevich's change of variables formula 
\cite{kontsevich}, \cite{denefloeser2}, \cite{denefloeser0}.  

\begin{cor}\label {change}
Let $\sigma: (M,E_0) \to (\R^d,0)$ be a real modification.  
Let $A\subset  \alphaanal (\R^d,0)$ be stable and suppose that 
$\ord jac (\sigma)$ is bounded  on $\sigma_*\inv (A)$. Then 
$$
\chic (A) = \int_{\sigma_*\inv (A)} (-1)^{- \ord jac (\sigma)} \, d\chic .
$$
\end{cor}

\begin{proof}
The function $\varphi = \ord jac (\sigma)$ is constructible on 
$A$.  Thus, by additivity of $\chic$, it suffices to show the formula 
only on $A_e:= A\cap \sigma_{*} (\Delta_e)$, for $e$ fixed.   
By Lemma \ref{keylemma}, for $n$ sufficiently large, $\sigma_{*n}$ is 
a locally trivial fibration over $\pi_n (A_e)$ with fiber
$\R^{e}$.  Hence $\chic (\sigma_{*n}\inv (\pi_n (A_e))) = 
 \chic (\R )^{-e} \chic  (\pi_n (A_e))$.  This ends the proof.  
\end{proof}

%We next show \eqref{} and \eqref{}.  

\medskip
\begin{proof}[Proof of \eqref{DL}, \eqref{DLpm}]

We show only \eqref{DL}.  The proof of \eqref{DLpm} is similar. 

The set  $\zed_n (f) = \pi_n\inv (\ixn (f) )$ is subanalytic and
$n$-stable.  The zeta function of $f$ can be  equivalently written as 
$$
\zetafa (T) = (-1)^d \sum_{n\ge 1} \chic (\zed_n(f)) T^n.
$$

Let $\zed_n (f\circ \sigma)= \sigma_* \inv (\zed_n (f))$
and $\zed_{n,e}(f \circ \sigma) 
= \zed_n(f \circ \sigma) \cap \Delta_e$.
Then $\zed_n(f \circ \sigma)$ is the disjoint union
of a finite number of  $\zed_{n,e}(f \circ \sigma)$.
Indeed, by comparing the multiplicities of $f \circ \sigma$ 
and $\jac \sigma$ along the components of the exceptional divisor 
we see that $\ord \jac \sigma \le n\max_i (\nu_i - 1)/N_i$ 
on $\zed_n (f\circ \sigma)$.   
(Here we use the assumption that the critical locus of $\sigma$
is contained in the zero set of $f$.
Otherwise the union may be infinite.)
By Kontsevich's change of variables formula
\begin{equation}\label{zed_f}
Z_f (T) = (-1)^d \sum_{n \ge 1} 
\sum_{e \le nq} (-1)^{-e}\chic (\zed_{n,e}(f \circ \sigma)) T^n
\end{equation}
where $q = \max_i (\nu_i - 1)/N_i$.

Fix $p\in  \Eo _I$.  
In a local system of coordinates at $p$ the germ of $f\circ \sigma$ at 
$p$, that we denote by $g: (\R^d,0)\to (\R,0)$,   is a normal 
crossings $g(y) = unit \cdot \prod_{i=1}^s y_i^{N_i}$, $s=|I|$.  
Let $\jac \sigma(y) = unit \cdot \prod_{i=1}^s y_i^{\nu_i -1}$.
We shall compute the weighted zeta function of $g$ that is 
$$
\hat Z_g (T) = (-1)^d \sum_{n \ge 1} 
\sum_{e \le nq} (-1)^{-e}\chic (\zed_{n,e}(g)) T^n, 
$$
where $\zed_{n,e}(g) = \zed_{n}(g)\cap \Delta_e$.  
Note that $\zed_{n,e}(g)$ is non-empty iff there are 
$k_1,\ldots, k_s$ such that $n=\sum k_i N_i$ and
$e=\sum k_i(\nu_i -1)$.  
We denote the set of such $k=(k_1,\ldots, k_s)$ by $A(n,e)$.  
Thus  $\zed_{n,e}(g)$ is the disjoin union 
$$
\zed_{n,e}(g) = \bigsqcup_{k\in A(n,e)} 
(\prod_{i=1}^s \zed_{k_i}(y_i^{N_i})
\times (\mathcal L (\R,0))^{d-s},
$$ 
and the last factor comes from the remaining $d-s$ variables $y_i$
that do not contribute to the zero of $g$.  Hence 
$$
\chic(\zed_{n,e}(g)) 
= (-1)^{d-s} \sum_{k\in A(n,e)} 
(\prod_{i=1}^s \chic (\zed_{k_i}(y_i^{N_i})).  
$$
Thus, by \eqref{ncrossing},   
\begin{eqnarray}\label{hatzeta}  
\hat Z_g (T) & = & (-1)^{s} \sum_{(k_1,\ldots, k_s)\in \N^s} 
\prod_{i=1}^{s} \chic (\zed_{k_i}(y_i^{N_i}))  
((-1)^{\nu_i -1}T^{N_i})^{k_i} \\
& = & \prod_{i=1}^{s}  \bigl ( \sum_k (-1)\chic (\zed_{k}(y_i^{N_i}))  
((-1)^{\nu_i -1}T^{N_i})^{k} \bigr )  \nonumber \\
& = &  (-2)^s \prod _{i=1}^s \frac {(-1)^{\nu_i }T^{N_i}} 
{1-(-1)^{\nu_i }T^{N_i}} \nonumber
\end{eqnarray}

Formula $\eqref{DL}$ follows now from $\eqref{zed_f}$
by integration (with respect to $\chic$) of \eqref{hatzeta} 
along the fibers of the projection
$\alphaanal(M,E_0) \to E_0 = \sigma^{-1}(0)$.
More precisely, to establish the equality of coefficients of $T^n$,
we integrate along the fibers of the projection
$\alphan (M,E_0) \to E_0$ restricted to
$\ixn (f \circ \sigma):= \sigma^{-1}_{*n} \ixn (f)$.
Then, by \eqref{hatzeta}, the Euler characteristic with compact
support of the fiber over $p \in \Eo_I$ is independent of the choice
of $p$ in $\Eo_I$. If we denote this Euler characteristic 
by $\chic (\ixn (f \circ \sigma)_I)$ then
$$
\chic (\ixn (f \circ \sigma)) = \sum_{I \neq \emptyset}
\chic (\Eo_I) \chic (\ixn (f \circ \sigma)_I),
$$
and the formula follows.
\end{proof}

\medskip
\begin{thm}\label{invariance}
Let $f, g :(\R^d,0) \to (\R, 0)$ be blow-analytically  equivalent 
function germs.  Then  $\zetafa = Z_g, 
\zetafp = Z_{g,+}, \zetafm = Z_{g,-}$.
\end{thm}

\begin{proof}
Since 
$f, g : (\R^d,0) \to (\R,0)$ are blow-analytically equivalent 
there are  real modifications
$\mu : (M,\mu^{-1}(0)) \to (\R^d,0)$,
$\mu^{\prime} : (M^{\prime},\mu^{\prime -1}(0)) \to (\R^d,0)$
and an analytic isomorphism $\Phi : (M,\mu^{-1}(0)) \to
(M^{\prime},\mu^{\prime -1}(0))$ such that 
 $f \circ \mu = g \circ\mu' \circ \Phi$.

First we show that we may assume that both $\mu$ and $\mu'$ satisfy
the properties required by Denef \& Loeser's formula.  
Let $\jac \mu$, $\jac (\mu' \circ \Phi)$  denote the jacobian
determinant of $\mu$, resp. of $\mu' \circ \Phi$.  
By \cite{hironaka}, \cite{bierstonemilman} there 
is a modification  $\mu_1: M_1\to M$ so that  
$f\circ \mu \circ \mu_1$ , $\jac \mu \circ \mu_1$, and 
$\jac (\mu' \circ \Phi) \circ \mu_1$  are normal crossings
simultaneously.  Moreover, we may assume that 
$\mu_1$ is a composition
of blowings-up with smooth centers that are in normal crossings
with the old exceptional divisors and hence that  $\jac \mu_1$
is normal crossings. Let $\sigma:= \mu \circ \mu_1$.  Then 
 $\jac \sigma (x) = \jac \mu_1 (x) \jac \mu (\mu_1(x))$ is normal 
crossings with respect to the same set of divisors.  
Set $\sigma' = \mu'\circ \Phi \circ \mu_1$.  Then $g\circ \sigma' = 
f\circ \sigma$ is normal crossings and so is $\jac \sigma' (x) = 
\jac \mu_1 (x) \jac (\mu'\circ \Phi) (\mu_1(x))$.  
Thus both $\sigma$ and $\sigma'$ satisfy the required properties.  

Let $E_i$ be an irreducible component of $(f \circ \sigma)^{-1}(0)$
(in $\sigma^{-1}(B_{\epsilon})$).
Since $g\circ \sigma' = f\circ \sigma$ the multiplicities
of these two functions coincide on $E_i$.  
Thus, by formulae \eqref{DL}, \eqref{DLpm}, 
 in order to show that the zeta 
functions of $f$ and $g$ 
%co\"{\i}ncide 
coincide it suffices to show that 
$\mult_{E_i} \jac \sigma$ and $\mult_{E_i} \jac \sigma'$ are 
of the same parity for any irreducible component $E_i$
of the exceptional divisor of $\sigma$
since $\mult \jac \sigma = 0$ outside the exceptional set 
$E$ of $\sigma$. 
Recall that $\Phi$ induces a homeomorphism 
$\phi : (\R^d,0) \to (\R^d,0)$
such that $f = g \circ \phi$. 
Then $\sigma (E)$ is of dimension $\le d-2$ and  
$\phi$ is analytic on the complement of 
$\sigma (E)$.  In particular the jacobian  $\jac \phi$ has constant 
sign on the complement of $\sigma (E)$.  
Fix $p\in  \Eo _i$ and a local system of coordinates on $M_1$ at $p$.  
Then, since $\phi \circ \sigma = \sigma'$, $\jac \phi (\sigma (x)) 
 \jac \sigma (x)= \jac \sigma' (x)$.  
In particular, 
$\jac \sigma (x)$ changes sign across $ \Eo _i$ iff so does 
$\jac \sigma' (x)$.  This shows that the multiplicities 
$\mult_{E_i} \jac \sigma$ and $\mult_{E_i} \jac \sigma'$ are of 
 the same parity, as claimed.  This ends the proof.  
\end{proof}

It follows that the modified zeta
functions $f$ and $g$ are also equal if $f$ and $g$ are
blow-analytically equivalent.

%%%%%%%%%%%%%%%%%%%%%%%%%%%%%%%%%%%%%%%%%%%%%%%%%%%%%%%%%%%%%%%%%

\bigskip
\section{Various Formulae to compute the Fukui Invariant}
\label{Fukui}
\medskip

\subsection{Formulae in terms of the resolution}\label{FoTRe}

Let $f : (\R^d,0) \to (\R,0)$ be an analytic function germ.
Take any analytic arc $\gamma : (\R,0) \to (\R^d,0)$.
Then $f(\gamma (t))$ is a convergent power series in $t$.
We denote by $\ord (f(\gamma (t)))$ its order in $t$.   
Set
$$
A(f) = \{ \ord (f(\gamma (t))) \in \N \cup \{ \infty \} 
; \gamma : (\R,0) \to (\R^d,0) \ C^{\omega} \} .
$$ 
In \cite{fukui}, T. Fukui proved that $A(f)$ is 
a blow-analytic invariant.
Namely, if analytic functions $f, g : (\R^d,0) \to (\R,0)$
are blow-analytically equivalent, then $A(f) = A(g)$.
We call $A(f)$ {\em the Fukui invariant}.
Note that the smallest number in $A(f)$ is the multiplicity of $f$.
For a positive integer $a\in \N$, set
$\N_{\ge a} = \{n\in \N; n\ge a\}$.

\begin{example}\label{example5}
Let $f : (\R^2,0) \to (\R,0)$ be a polynomial function defined
by $f(x,y) = x^3 - y^5$.
Then
$$
A(f) = 3\N \cup 5\N \cup \N_{\ge 16} \cup \{ \infty \}
= \{ 3, 5, 6, 9, 10, 12, 15, 16, 17, \cdots \} \cup \{ \infty \}.
$$
Any integer $15+s \in A(f)$, $s \in \N$, is attained 
along  $\gamma (t) = (t^5+t^{5+s},t^3)$.
\end{example}

For an analytic function germ $f : (\R^d,0) \to (\R,0)$,
let $\sigma: M\to \R^d$ be a simplification of $f^{-1}(0)$,
namely, $\sigma$ is a composition of a finite number of blowings-up,
$M$ is smooth and $f\circ \sigma$ is normal crossing.
As in Subsection 1.2, we denote by $E_i$, $i\in J$, 
the irreducible components of $(f\circ \sigma) \inv (0)$ 
(in $\sigma \inv (B_\varepsilon)$, where $B_\varepsilon$ 
is a small ball in $\R^d$ centered at the origin).  
%We suppose also that $\sigma \inv (0)$ is 
%the union of some of $E_i$.   
For each $i\in J$, let $N_i= \mult _{E_i} f\circ \sigma$.  
Denote for $I\subset J$,  
$E_I = \bigcap _{i\in I} E_i$ 
and $\Eo _I = E_I \setminus \bigcup _{j\in J\setminus I} E_j$.
We put
$$
\mathcal{C} = \{ I ; \Eo _I \cap \sigma^{-1}(0)
\ne \emptyset \}.
$$

\begin{rem}\label{remark0}
As stated in Section 1, we can assume that 
$\sigma \inv (0)$ is the union of some of $E_i$.
Then $\mathcal{C} = \{ I | E_I \subset \sigma^{-1}(0) \}$.
\end{rem}

For $A$, $B \subset \N \cup \{ \infty \}$,
define $A + B = \{ a+b \in \N \cup \{ \infty \};
a \in A, b \in B \}$,
where we set $a+b = \infty$ if $a = \infty$
or $b = \infty$. Let us put 
$$
\Omega_I(f) = (N_{i_1}\N + \cdots + N_{i_p}\N) \cup \{ \infty \},
$$
for $I = (i_1, \cdots , i_p) \in \mathcal{C}$. 
%Then we can write down the Fukui invariant $A(f)$,
%using $\Omega_I(f)$.

\medskip
\begin{thm}\label{IzuKoKu1} (\cite{izumikoikekuo})
Let $f : (\R^d,0) \to (\R, 0)$ 
be an analytic function germ and let
$\sigma$ be a simplification of $f^{-1}(0)$. 
Then we have
$$
A(f)=\bigcup_{I \in \mathcal{C}} \Omega_I(f).
$$
\end{thm}

Let us put
\begin{eqnarray*}
& \mathcal{C}^+ := \{ I \in {\mathcal C} ; \Eo _I \cap 
\sigma^{-1}(0) \cap \overline{P(f)} \neq \emptyset \} , \quad 
& P(f) :=  \{ x \in M ; f \circ \sigma(x) > 0 \}, \\  
& \mathcal{C}^- :=  \{ I \in {\mathcal C} ; \Eo _I \cap 
\sigma^{-1}(0) \cap \overline{N(f)} \neq \emptyset \} , \quad 
& N(f) := \{ x \in M ; f \circ \sigma(x) < 0\}, 
\end{eqnarray*}
where the overlines denote the closures in $M$. 

Let $\lambda : U \to \R^d$ be an analytic arc with $\lambda (0) = 0$, 
where $U$ denotes a neighborhood of $0 \in \R$. 
We call $\lambda$  {\it nonnegative} (resp. {\it nonpositive}) 
{\it for} $f$ if 
$(f \circ \lambda)(t) \geq 0$ (resp. $\leq 0$) 
in a positive half neighborhood $[0,\delta) \subset U$. 
Then we define {\em the Fukui invariants with sign} by 

\vspace{3mm}

\qquad $A_+(f) := \ \{ \ord (f \circ \lambda) ; \lambda$
is a nonnegative arc through $0$ for $f \}$, 

\qquad $A_-(f) := \ \{ \ord (f \circ \lambda) ; \lambda$ 
is a nonpositive arc through $0$ for $f \}$,

\vspace{3mm}

\noindent respectively. It is easy to see that
these $A_+(f)$ and $A_-(f)$ are also blow-analytic invariants.
Remark that $A(f)=A_+(f) \cup A_-(f)$. 
Then we have the following formulae to compute
the Fukui invariants with sign:

\medskip
\begin{thm}\label{IzuKoKu2} (\cite{izumikoikekuo})
Let $f : (\R^d,0) \to (\R,0)$ be an analytic function germ.
Then we have
$$
A_+(f) = \bigcup_{I \in {\mathcal C}^+} \Omega_I(f), 
\  A_-(f) = \bigcup_{I \in {\mathcal C}^-} \Omega_I(f).
$$
\end{thm}

\medskip

\subsection{List of the Fukui invariants for $\pm x^p \pm y^q$}\label{List}

Let $p, q \in \N$, and let $(p,q) = d$.
Here, $(p,q)$ denotes $\gcd(p,q)$.
Then there are $p_1, q_1 \in \N$
such that $p = p_1d$, $q = q_1d$ and $(p_1,q_1) = 1$.
Set $[p,q] = LCM(p,q) = p_1q_1d = pq_1 = p_1q$.

Using the argument of example \ref{example5}, we compute  
the Fukui invariants for Brieskorn polynomials $f(x,y) = \pm x^p \pm y^q$, 
$(x,y)\in \R^2$, $p \le q$, listed in the table below. 

\begin{table}[hbt]
\begin{tabular}{|l|l|}
\hline 
$f(x,y)$ & Fukui invariants \\
\hline \hline 
$\pm x ^{p} \pm y^{q}$, \, p, q odd &  $A(f) = A_+(f) = A_-(f) = 
p\N \cup q\N \cup 
\N_{\ge [p,q]} \cup \{ \infty \}$    \\ 
% $\pm x ^{p} \pm y^{q} $ & $A_+(f) = A_-(f) = A(f)$ \hfill\\
\hline \hline
 p odd, q even  \quad & $A(f) = p\N \cup q\N \cup \N_{\ge [p,q]} \cup 
\{ \infty \}$  \\ \hline
 $\pm x ^{p} + y^{q}$ & $A_+(f) =A(f)$, $A_-(f) = p\N  \cup \N_{\ge
 [p,q]} \cup 
\{ \infty \} \qquad $ \\ \hline
$\pm x ^{p} - y^{q}$ & $A_-(f) = A(f)$, $A_+(f) = p\N  \cup \N_{\ge
 [p,q]} \cup \{ \infty \}$ \\ \hline \hline
 p even, q odd &  $A(f) = p\N \cup q\N \cup \N_{\ge [p,q]}  
\cup \{ \infty \}$   \\  \hline
 $ x ^{p} \pm y^{q}$ & $A_+(f) =A(f)$, $A_-(f) = q\N \cup \N_{\ge
 [p,q]} \cup \{ \infty \}$   \\ \hline
 $- x ^{p} \pm y^{q}$ & $A_-(f) = A(f)$, $A_+(f) = q\N \cup \N_{\ge [p,q]} 
\cup \{ \infty \}$ \\ \hline \hline
$\pm(x ^{p} - y^{q})$, \, p,q even \, &  $A(f) = p\N \cup q\N \cup
 \N_{\ge [p,q]}
\cup \{ \infty \}$   \\  
% $\pm(x ^{p} - y^{q})$ &   \\ 
\hline
$ x ^{p} - y^{q}$ & $A_+(f) = p\N \cup \N_{\ge [p,q]} \cup \{ \infty \}$, 
\, $A_-(f) = q\N \cup \N_{\ge [p,q]} \cup \{ \infty \}$  \,  
% \\ & $A_-(f) = q\N \cup \N_{\ge [p,q]} \cup \{ \infty \}$   
\\ \hline
$- x ^{p} + y^{q}$ & $A_+(f) = q\N \cup \N_{\ge [p,q]} \cup \{ \infty \}$, 
\, $A_-(f) = p\N \cup \N_{\ge [p,q]} \cup \{ \infty \}$ \, 
% \\ & $A_-(f) = p\N \cup \N_{\ge [p,q]} \cup \{ \infty \}$   
\\ \hline \hline
$\pm(x ^{p} + y^{q})$, \, p,q even \, &  $A(f) = p\N \cup q\N  \cup 
\{ \infty \}$   \\  
%$\pm(x ^{p} + y^{q})$ &   \\ 
\hline
$ x ^{p} + y^{q}$ & $A_+(f) = A(f)$, $A_-(f) =\{ \infty \}$    
\\ \hline
$- x ^{p} - y^{q}$ & $A_-(f) = A(f)$, $A_+(f) = \{ \infty \}$ \\
\hline \hline
\end{tabular}
%\caption{}
\end{table}

\begin{rem}\label{remark1}
Let $f_1(x,y) = \pm x^p + y^q$ and $f_2(x,y) = \pm x^p - y^q$, 
$p$ odd, $q$ even.
If $q$ is divisible by $p$, then $[p,q] = q = q_1p$.
Thus $A(f_1) = A_{\pm}(f_1) = A(f_2) = A_{\pm}(f_2)$. 

If $q$ is not divisible by $p$, then $[p,q] > q$.
Thus $A_+(f_1) \ne A_+(f_2)$ and $A_-(f_1) \ne A_-(f_2)$.
\end{rem}

\subsection{Thom-Sebastiani formulae for the Fukui invariant}
\label{TSFukui}

Let $f : (\R^{d_1},0) \to (\R,0)$ and $g : (\R^{d_2},0) \to 
(\R,0)$ be analytic function germs. Define
$f*g : (\R^{d_1+d_2},0) \to (\R,0)$ by
$(f*g)(x,y) := f(x) +g(y)$ as in Section 2, 
and define also $f \cdot g : (\R^{d_1+d_2},0) \to (\R,0)$
by $(f \cdot g)(x,y) := f(x) \times g(y)$.
In this subsection, we give the Thom-Sebastiani formulae
expressing the Fukui invariants of $f(x) + g(y)$ and
$f(x) \times g(y)$ in terms of the Fukui invariants of
$f(x)$ and $g(y)$.

\medskip
\begin{thm}\label{Fukui*}
Let $M_1 = \min(A_+(f) \cap A_-(g))$ and
$M_2 = \min(A_-(f) \cap A_+(g))$.
\begin{eqnarray}\label{A(f*g)}
& A(f*g) = & A(f) \cup A(g) \cup (M_1+\N) \cup (M_2+\N), \\
\label{A_+(f*g)}
& A_+(f*g) = & A_+(f) \cup A_+(g) \cup (M_1+\N) \cup (M_2+\N), \\
\label{A_-(f*g)}
& A_-(f*g) = & A_-(f) \cup A_-(g) \cup (M_1+\N) \cup (M_2+\N).
\end{eqnarray}
\end{thm}

\begin{proof}
We show only \eqref{A(f*g)}.

\vspace{2mm}

\noindent($\subset$) Take any $k \in A(f*g)$.
We may assume that $k < \infty$ since
$\infty \in A(f)$ or $A(g)$.
Then there is an analytic arc $\nu = (\lambda,\mu) :
(\R,0) \to (\R^{d_1} \times \R^{d_2},(0,0))$
such that $\ord ((f*g) \circ \nu) = k$.   Let
\begin{eqnarray*}
& (f \circ \lambda)(t) = & a_u t^u + a_{u+1} t^{u+1} + \cdots , 
\ a_u \ne 0,  \\
& (g \circ \mu)(t) = & b_v t^v + b_{v+1} t^{v+1} + \cdots ,
\ b_v \ne 0. 
\end{eqnarray*}
Then $u = \ord (f \circ \lambda) \in A(f)$ and
$v = \ord (g \circ \mu) \in A(g)$.
Since $k \in A(f*g)$ and, $u \le k$ or $v \le k$, 
it suffices to consider the following three cases:

(i) $u = k$ and $v \ge k$; In this case, $k \in A(f)$.

(ii) $u \ge k$ and $v = k$; In this case, $k \in A(g)$.

(iii) $u = v < k$; In this case, $u \in A_+(f)$ and $v \in A_-(g)$,
or $u \in A_-(f)$ and $v \in A_+(g)$.
This means 
$$
u = v \in (A_+(f) \cap A_-(g)) \cup (A_-(f) \cap A_+(g)).
$$
It follows that $k > u = v \ge \min(M_1,M_2)$.

If $k \le \min(M_1,M_2)$, then case (i) or case (ii) holds. Thus
$$
k \in A(f) \cup A(g) \cup (M_1+\N) \cup (M_2+\N)
$$
because $k > \min(M_1,M_2)$ implies
$k \in (M_1+\N) \cup (M_2+\N)$.

\vspace{2mm}

\noindent
($\supset$) It is obvious that $A(f), \ A(g) \subset A(f*g)$. 
Let us show $M_1 + \N \subset A(f*g)$. 

First we recall the reparametrization formulae of remark 
\ref{trivial1} and lemma \ref{trivial}.
Let $h : (\R,0) \to (\R,0)$ be an analytic function defined by
$$
h(t) = a_k t^k + a_{k+1} t^{k+1} + \cdots , \ a_k \ne 0.
$$
Then, if we replace $t$ by $\alpha t$, $\alpha \ne 0$,  
$$
h(\alpha t) = a_k \alpha^k t^k + a_{k+1}^{\prime} t^{k+1} 
+ \cdots .
$$
Let $A = \{a_k \alpha^k; \alpha \in \R^* \}$. Then $A = \R^*$ for 
$k$ odd, $A = \R_{> 0}$ for $k$ even and $a_k > 0$,
and $A = \R_{< 0}$ for $k$  even and $a_k <0$. 
Similarly, if we replace $t$ by $t+\alpha t^{i+1}$, $i\ge 1$, 
$$
h(t+\alpha t^{i+1}) = a_k t^k + \cdots + a_{k+i-1} t^{k+i-1} +
(a_k k \alpha + a_{k+i}) t^{k+i} %+ a_{k+i+1}^{\prime} t^{k+i+1} 
+ \cdots ,
$$
and in this case $\{a_k k \alpha + a_{k+i}; \alpha \in \R \} = \R$.

Take $k+j \in M_1 + \N$ such that $k = M_1$ and $j \in \N$.
Then there are a nonnegative arc for $f$, 
$\lambda : (\R,0) \to (\R^{d_1},0)$,
and a nonpositive arc for $g$,
$\mu : (\R,0) \to (\R^{d_2},0)$,
such that $\ord (f \circ \lambda) = \ord (g \circ \mu) = k$.
Then 
\begin{eqnarray*}
& (f \circ \lambda)(t) =  a_k t^k + a_{k+1} t^{k+1} + \cdots , \quad &a_k>0, 
\\
& (g \circ \mu)(t) =  b_k t^k + b_{k+1} t^{k+1} + \cdots , \quad & b_k <0. 
\end{eqnarray*}
By the above there is a reparametrization 
$\mu^{(1)} : (\R,0) \to (\R^{d_2},0)$ of $\mu$  such that
$$
(g \circ \mu^{(1)})(t) = - a_k t^k + b_{k+1}^{(1)} t^{k+1} + \cdots .
$$
Using the second type of reparametrizations we can construct 
by induction on $i$ an analytic arc $\mu^{(i)} : (\R,0) \to
(\R^{d_2},0)$ such that
$$
(g \circ \mu^{(i)})(t) = - a_k t^k - a_{k+1} t^{k+1}
- \cdots - a_{k+i-1} t^{k+i-1} + b_{k+i}^{(i)} t^{k+i} + \cdots ,
$$
for $2 \le i \le j$.
Using the same argument again, we show that there is an analytic arc
$\check{\mu} : (\R,0) \to (\R^{d_2},0)$ such that
$$
(g \circ \check{\mu})(t) = - a_k t^k - a_{k+1} t^{k+1}
- \cdots - a_{k+j-1} t^{k+j-1} + \check{b}_{k+j} t^{k+j} + \cdots 
$$
with $a_{k+j} + \check{b}_{k+j} \ne 0$.
Define $\nu : (\R,0) \to (\R^{d_1} \times \R^{d_2},(0,0))$
by $\nu(t) = (\lambda (t),\check{\mu}(t))$.
Then $\ord ((f*g) \circ \nu) = k+j$. 
Thus $k+j \in A(f*g)$, namely, $M_1 + \N \subset A(f*g)$.

We can similarly show $M_2 + \N \subset A(f*g)$.
\end{proof}

\begin{example}\label{example6}
Let $f(x) = x^4$ and $g(y) = y^6$.
Then $A(f) = A_+(f) = 4\N \cup \{ \infty \}$,
$A(g) = A_+(g) = 6\N \cup \{ \infty \}$,
$A_-(f) = A_-(g) = \{ \infty \}$ and
$M_1 = M_2 = \infty$.

Thus $A(f*g) = A_+(f*g) = 4\N \cup 6\N \cup \{ \infty \}$ and
$A_-(f*g) = \{ \infty \}$.
\end{example}

Concerning the Fukui invariant for $f \cdot g$, 
we can easily show following formulae.

\begin{prop}\label{Fukuimulti}
\begin{eqnarray}\label{A(f.g)}
& A(f \cdot g) = & A(f) + A(g), \\
\label{A_+(f.g)}
& A_+(f \cdot g) = & (A_+(f) + A_+(g)) \cup (A_-(f) + A_-(g)), \\
\label{A_-(f.g)}
& A_-(f \cdot g) = & (A_+(f) + A_-(g)) \cup (A_-(f) + A_+(g)).
\end{eqnarray}
\end{prop}

\begin{rem}\label{remark3}
%Let $f : (\R,0) \to (\R,0)$ be an analytic function.  Then 
$A(f) = (\min A(f))\N \cup \{ \infty \}$, 
$A_+(f) = (\min A_+(f))\N \cup \{ \infty \}$ and
$A_-(f) = (\min A_-(f))\N \cup \{ \infty \}$.
\end{rem}

\begin{example}\label{example7}
Let $f(x,y) = cx^py^q$, $c \ne 0$, $x,y \in \R^2$. Then
$$
A(f) = \{ ap + bq; a, b \in \N \} \cup \{ \infty \}.
$$
\noindent (1) Let $p$ or $q$ be odd. Then $A_+(f) = A_-(f) = A(f)$.

\noindent (2) Let $p$ and $q$ be even.

(i) If $c > 0$, $A_+(f) = A(f)$ and $A_-(f) = \{ \infty \}$.

(ii) If $c < 0$, $A_+(f) = \{ \infty \}$ and $A_-(f) = A(f).$

\end{example}

%\medskip

%%%%%%%%%%%%%%%%%%%%%%%%%%%%%%%%%%%%%%%%%%%%%%%%%%%%%%%%%%%%%

\bigskip
\section{Two variables Brieskorn polynomials}
\label{2BP}
\medskip

\subsection{Classification of two variables Brieskorn polynomials}
\label{classification2BP}

Let $f : (\R^2,0) \to (\R,0)$ be a two variables Brieskorn polynomial
defined by $f(x,y) = \pm x^p \pm y^q$, $p \le q$.
If $0 \in \R^2$ is a regular point of $f$, i.e.
$p = 1$, then f is analytically equivalent to 
$g(x,y) = x$ by the Implicit Function Theorem.
After this, we assume that $0 \in \R^2$ is a singular
point of $f$, i.e. $2 \le p \le q$.

Let $\N_e$ (resp. $\N_o$) denote the set of positive even integers
(resp. positive odd integers).
Set
\begin{eqnarray*}
& \mathfrak{M}: = & \{ (p,q) \in \N_{\ge 2} \times \N_{\ge 2} 
; p \le q \},  \\ 
& \mathfrak{N}: = & \mathfrak{M} - \{ (p,mp) \in \mathfrak{M} 
; p \in \N_o, \ m \in \N_e \}.
\end{eqnarray*}

Let us consider the classification of Brieskorn polynomials
by blow-analytic equivalence.
We denote by $(\pm x, \pm y)$ the Klein group $G = \Z_2 \oplus \Z_2$
consisting of the following four transformations of $\R^2$:
$$
(x,y) \to (x,y), \ (x,y) \to (-x,y), \ (x,y) \to (x,-y),
\ (x,y) \to (-x,-y).
$$
For a subset $A$ of $\{ f(x,y) = \pm x^p \pm y^q 
\ | \ (p,q) \in \mathfrak{M} \}$, let $A/b.a.e$ (resp.
$A/(\pm x,\pm y)$) denote
the quotient of $A$ by blow-analytic equivalence
(resp. the Klein $G$-equivalence).
Then we have the following blow-analytic classification.

\medskip
\begin{thm}\label{classification} \ \ \
$\{ f(x,y) = \pm x^p \pm y^q ; (p,q) \in \mathfrak{M} \}/b.a.e.$

\vspace{2mm}

\noindent $= \{ f(x,y) = \pm x^p \pm y^q ; (p,q) \in \mathfrak{N} 
\}/(\pm x,\pm y) \cup
\{ x^p + y^{mp} ; p \in \N_o \cap \N_{\ge 2}, \ m \in \N_e \}.$
\end{thm}

\begin{proof}
By our list of the Fukui invariant in Subsection \ref{List},
we can distinguish all real Brieskorn polynomials of two variables
$f(x,y) = \pm x^p \pm y^q$, $(p,q) \in \mathfrak{M}$,
up to $\{ (\pm x,\pm y) \}$ by the Fukui invariant
except the following two cases:

\vspace{2mm}

\noindent Case (i): $x^p + y^{mp}$ \ for a fixed even $p$ 
and $m = 1, 2, 3, \cdots$, 

or $- x^p - y^{mp}$ \ for a fixed even $p$
and $m = 1, 2, 3, \cdots$.

\noindent Case (ii): $\pm x^p + y^{mp}$ and $\pm x^p - y^{mp}$
\ for fixed odd $p \ge 3$ and even $m$.

\vspace{2mm}

We first consider case (i).
For a fixed even $p$, let $f_m(x,y) = x^p + y^{mp}$
and $g_m(y) = y^{mp}$,
$m = 1, 2, 3, \cdots $.
In this case, 
$$
A(f_m) = A_+(f_m) = \{ p, 2p, 3p, \cdots \} \cup \{ \infty \},
\ \ A_-(f_m) = \{ \infty \}, \ \ m = 1, 2, 3, \cdots.
$$
Since $p$ and $mp$ are even,
it follows from corollary \ref{xmeven} that 
if $Z_{x^p*g_m}(T) = Z_{x^p*g_n}(T)$,
then $Z_{g_m}(T) = Z_{g_n}(T)$.
On the other hand, as seen in example 1.3.1, 
$Z_{g_m}(T) \neq Z_{g_n}(T)$ if $m \ne n$.
Since the zeta function is a blow-analytic invariant,
$f_m$ and $f_n$ are not blow-analytically equivalent if $m \ne n$.
The case of $- x^p  - y^{mp}$ follows similarly.

We next consider case (ii). In this case,
$x^p + y^{mp}$ (resp. $x^p - y^{mp}$) is equivalent to
$- x^p + y^{mp}$ (resp. $- x^p - y^{mp}$)
under the transformation of $\R^2$:
$(x,y) \to (-x,y)$.
Therefore, we treat only $f(x,y) = x^p + y^{mp}$
and $g(x,y) = x^p - y^{mp}$ 
for fixed odd $p \ge 3$ and even $m$.
Remark that the Fukui invariants $A(f)$ and $A_{\pm}(f)$
and the zeta functions $Z_f(T)$ and $Z_{f,\pm}(T)$ 
coincide with $A(g)$, $A_{\pm}(g)$,
$Z_g(T)$ and $Z_{g,\pm}(T)$, respectively. 

Here we recall the Fukui-Paunescu Theorem.

\begin{lem}\label{FuPa}
(T. Fukui - L. Paunescu \cite{fukuipaunescu},
T. Fukui - E. Yoshinaga \cite{fukuiyoshinaga})
Given a system of weights $\alpha = (\alpha_1,\cdots,\alpha_d)$.
Let $f_s : (\R^d,0) \to (\R,0), \ s \in I = [0,1],$ be an analytic 
family of analytic function germs.
Suppose that for each $s \in I$, the weighted initial form
of $f_s$ with respect to $\alpha$ is of the same weighted degree
and has an isolated singularity at $0 \in \R^d$.
Then $\{ f_s \}_{s \in I}$ is blow-analytically trivial
over $I$.
\end{lem}

Let $\{ f_s \}$ be a family of polynomial functions defined by
$$
f_s(x,y) = x^p + pxy^{m(p-1)} + sy^{mp}, \ \ s \in [-1,1].
$$
Then it follows from lemma \ref{FuPa} that
$x^p + pxy^{m(p-1)} + y^{mp}$ and 
$x^p + pxy^{m(p-1)} - y^{mp}$ are blow-analytically equivalent.

Nextly, let $\{ g_s \}$ and $\{ h_s \}$ be families of polynomial
functions defined by
\begin{eqnarray*}
g_s(x,y) = x^p + psxy^{m(p-1)} + y^{mp}, \ \ s \in [0,1], \\
h_s(x,y) = x^p + psxy^{m(p-1)} - y^{mp}, \ \ s \in [0,1].
\end{eqnarray*}
Then, by the same reason as above,
$x^p + pxy^{m(p-1)} + y^{mp}$ (resp. $x^p + pxy^{m(p-1)} - y^{mp}$)
are blow-analytically equivalent to $x^p + y^{mp}$
(resp. $x^p - y^{mp}$).
Since blow-analytic equivalence is an equivalence relation (\cite{kuo2}),
$x^p + y^{mp}$ and $x^p - y^{mp}$ are blow-analytically equivalent.

This completes the proof of the theorem.
\end{proof}

Concerning cases (i) and (ii) in the proof of Theorem 6.1,
we have the following remarks.

\begin{rem}\label{remark4}
By the above proof, we see that the Fukui invariants
distinguish all real Brieskorn polynomials of two variables
except case (i) and are not enough to give a complete classification 
of Brieskorn polynomials by blow-analytic equivalence.
Then it gives rise to the following natural question:

\vspace{3mm}

{\em Is the blow-analytic type of Brieskorn polynomials 
completely determined by the zeta functions?}

\vspace{3mm}

The answer is \lq No'.
Our zeta functions distinguish the blow-analytic types of
all real Brieskorn polynomials of two variables
except $f_m^+(x,y) = x^{2m} + y^{2m}$
and $f_m^-(x,y) = - (x^{2m} + y^{2m})$,
$m = 1, 2, 3, \cdots$.
We shall see this fact in the next subsection.
As seen in example \ref{2k2k},
$Z_{f_m^{\pm}}(T) = Z_{f_m^{\pm},\pm}(T) = 0$
for any $m$.
On the other hand, $A(f_m^{\pm})$ is different
from $A(f_n^{\pm})$ if $m \ne n$
and $A_+(f_m^+) \neq A_+(f_m^-)$.

These remarks mean that the Fukui invariants and the zeta functions
are compensating each other for our blow-analytic classification. 
\end{rem}

\begin{rem}\label{remark5}
Consider two functions of case (ii),
$f(x,y) = x^p + y^{mp}$ and $g(x,y) = x^p - y^{mp}$
\ for fixed odd $p \ge 3$ and even $m$.
These functions are exceptional in our classification
since they are blow-analytically equivalent,
but not Klein $G$-equivalent.
It is easy to see that they are not
analytically equivalent.
In addition, it was shown recently in
\cite{henryparusinski2}, \cite{henryparusinski1}
that $f$ and $g$ are not even bi-Lipschitz equivalent.
\end{rem}

\subsection{Distinction of Brieskorn polynomials by zeta functions.}
\label{distinction}

Let $f(x,y) = \pm x^p \pm y^q$, $ 2 \le p \le q$.
Considering $f(x,y)$ up to Klein $G$-equivalence,
we assume the following:

\vspace{2mm}

\noindent (i) In case $p$ (resp. $q$) is odd, we consider
only the positive case that is the coefficient at $x^p$
(resp. $y^q$) is $+1$.

\noindent (ii) In case $p = q$ are even, we consider 
$f(x,y) = x^p - y^q$ but not $f(x,y) = - x^p + y^q$.

\vspace{2mm}

We show that our zeta functions distinguish all real 
Brieskorn polynomials of two variables up to blow-analytic
equivalence except
\begin{equation}\label{sumsquare}
f(x,y) = \pm (x^p + y^p), \ \ p = 2, 4, 6, \cdots .
\end{equation}
Note that $Z_f(T) = Z_{f,\pm}(T) \equiv 0$ only for Brieskorn polynomials
of form (\ref{sumsquare}) in the two variables case.

Assume that $f(x,y) = \pm x^p \pm y^q$ is not of form 
(\ref{sumsquare}). Let $Z_f(T) = \sum_{i \ge 1} c_i T^i$,
$Z_{f,\pm}(T) = \sum_{i \ge 1} c_i^{\pm} T^i$ as above.
Then, by theorem \ref{Thom-Seb} and example \ref{xtom}, 
we see that $c_i = 0, \ 1 \le i \le p - 1$, and $c_p \neq 0$.
Therefore $p$ is determined by $Z_f(T)$.

We first consider the even case that is $p$ is even.
If $f(x,y) = x^p \pm y^q$ (resp. $- x^p \pm y^q$),
$p < q$, then $c_p^+ = c_p \neq 0$, $c_p^- = 0$
(resp. $c_p^- = c_p \neq 0$, $c_p^+ = 0$).
Therefore the sign at $x^p$ is determined by $Z_{f,\pm}(T)$.
Let $\phi (x) = \pm x^p$.
By corollary \ref{xmeven}, $Z_{y^q,\pm}(T)$
(resp. $Z_{-y^q,\pm}(T)$) or $\mzeta_{y^q,\pm}(T)$
(resp. $\mzeta_{-y^q,\pm}(T)$) can be computed from
$Z_{\phi *y^q,\pm}(T)$ 
(resp. $Z_{\phi *(-y^q),\pm}(T)$).
As seen in example \ref{example2},
$\mzeta_{\pm y^q,\pm}(T)$ are different from
$\mzeta_{\pm y^{q'},\pm}(T)$ if $q \neq q'$, 
and $\mzeta_{y^q,\pm}(T)$ are different from
$\mzeta_{-y^q,\pm}(T)$ (if $q$ is even).
Therefore $q$ and the sign at $y^q$ are determined by $Z_{f,\pm}(T)$.

We next consider the odd case. Then, by proposition \ref{xmodd},
$q$ is determined by $\mzeta_{f,\pm}(T)$.
If $q$ is even and not divisible by $p$,
the sign at $y^q$ is also determined by $\mzeta_{f,\pm}(T)$.
On the other hand, as shown in the preceding subsection,
if $q$ is even and divisible by $p$,
$x^p + y^q$ and $x^p - y^q$ are blow-analytically equivalent.
Therefore the zeta functions distinguish Brieskorn polynomials
up to blow-analytic equivalence in this case, too.

%%%%%%%%%%%%%%%%%%%%%%%%%%%%%%%%%%%%%%%%%%%%%%%%%%%%%%%%%%%%%%%%%%

\bigskip
\section{Examples in three variables}
\label{3V}
\medskip

\subsection{Brieskorn polynomials of three variables}
\label{3B}

Using the zeta functions and the Fukui invariants we classify 
blow-analytic types of Brieskorn polynomials of three variables, except 
for the following families: 
$\{x^p+y^{kp} +z^{kp}; k\in \N \}$, $\{-(x^p+y^{kp} +z^{kp}); k\in \N \}$, 
$p$ even.  

The following proposition generalizes proposition \ref{xmodd}.  

\begin{prop}\label{brieskorn}
Let $f(x_1,\ldots,x_d)$ be a Brieskorn polynomial, $f(x_1,\ldots,x_d)
= \pm x_1^{m_1} \pm \cdots \pm x_d^{m_d}$, all $m_i \ge 2$, 
and let $g(y) = \pm\ y^{r}$.  
Then $r$ is determined by the zeta functions of $f$ and of
$f*g$.  If, moreover, $r$ is even and $r\notin \bigcup_{m_i \text
{odd}} m_i \N$ then the sign at $y^r$ is determined, too.
\end{prop}

\begin{proof}
We use notation \eqref{mzeta} for the modified zeta functions of 
$f$, $g$, and $f*g$.  By assumption the coefficients $\tilde A^\pm _n$, 
resp. $\tilde C^\pm _n$, of the modified zeta functions of $f$, resp. 
$f*g$ are given.  Hence, by Thom-Sebastiani Formulae \eqref{tildeC}, 
we may determine the
coefficients $\tilde B^\pm_n$ of the modified zeta functions of
$g$ for all $n$ such that $\tilde C^\pm_n =0$ that is for
$n\in U:= \N \setminus \bigcup_{m_i \text{odd}} m_i \N$.

If there is $n\in U$ such that $\tilde B^+_n =0$ then $r$ is 
odd and equals the minimum of such $n$.  Similarly, if
there is $n\in U$ such that $\tilde B^+_n \ne \tilde B^-_n$
then,  $r$ is even and equals the minimum of such $n$.  In this
case we may determine the sign in $g(y) = \pm y^{r}$.

From now on we suppose that 
$$
\tilde B^+_n = \tilde B^-_n \ne 0 \qquad \text{for all } n\in U. 
$$
Then $r$ is a multiple of one of odd $m_i$'s.  We shall show that
the values 
\begin{equation}\label{Bn}
\tilde B^\pm_n, \, n \in U
\end{equation}
 determine $r$.  Without loss
of generality we may suppose that all odd $m_i$ are distincts
prime numbers.  Otherwise, without increasing $U$, we replace the
set of odd $m_i$'s by the set of all their prime divisors.  
Thus we assume $U = \N \setminus \bigcup_{p\in P} p\N$, where 
$P$ is a finite set of odd prime numbers.  Let $m$ be the product of 
all $p\in P$.  

First we show that $m^{\prime}=(m,r)$ is determined by the coefficients 
\eqref{Bn}.  
Let $m=m' m''$.  Then $(m'',r)=1$.  So there exist $a, b\in \Z$ such 
that for all $k\in \N$  
$$
(a+km'')r = (kr -b)m'' -1.
$$
Since $m''$ is odd, choosing $k$ we may suppose that $a+km''$ is even, and 
$a+km'', kr -b \in \N$ if $k$ is sufficiently large.  Then $kr -b$ is 
odd.  Fix such natural 
$$
q= Ar = Bm''-1 \qquad A \text{ even }, B \text { odd}.  
$$
Each $p\in P$ divides either $r$ or $m''$ and hence does not divide 
$q-1$ nor $q+2$, i.e. $q-1, q+2 \in U$.  Thus, 
by example \ref{example2}, if $g(y) = \pm y^r$ then 
\begin{equation}\label{q-1q+2}
\tilde B^+_{q-1}=-1, \,  \tilde B^+_{q+2}=1.
\end{equation}
Suppose now that $g(y) = \pm y^{r_1}$ gives the same coefficients \eqref{Bn} 
as $g(y) = \pm y^{r}$ and that there is $p_0 \in P$ 
such that $p_0$ divides $r_1$ but it does not divide $r$.  We show that this 
contradicts \eqref{q-1q+2}.  Note that \eqref{q-1q+2} is possible 
only if either $q$ or $q+1$ is an even multiple of $r_1$.  
%(neither $q-1$ nor $q+2$ is divisible by $r_1$ so it is not divisible 
%by $p_0$).  
Firstly, $q+1=Bm''$, as an odd number, cannot be an even 
multiple of 
$r_1$.  Secondly, $p_0$ divides $q+1=Bm''$ so it does not divide $q$.
Hence $r_1$ cannot divide $q$.  Thus if $g(y) = \pm y^{r}$ and 
$g(y) = \pm y^{r_1}$ 
give the same coefficients \eqref{Bn} they have the same factors in $P$.  That 
is $(m,r) = (m,r_1)$.  

Let  $m'=(m,r) = (m,r_1)$, $m=m'm''$.  Then $(r,r_1) = dm'$ where $(d,m'')=1$.  Suppose $r\ne r_1$. Then one of them, say $r_1$, is strictly 
bigger than $dm'$.  By assumptions $(r_1,rm'')=dm'$ so there is $q \in \N$ of the form 
$$
q = Arm'' = Br_1 +dm'.
$$
Clearly $q$ is a multiple of $m$ so $q-1,q+1\in U$.  But then, for $g=\pm y^r$, $$
\tilde B^+_{q-1} =- \tilde B^+_{q+1}\ne 0. 
$$
But this is not possible for $g=\pm y^{r_1}$ 
since $Br_1 <q-1 <q+1 <(B+1)r_1$.  
This ends the proof.
\end{proof}

%Now we make one remark on the Fukui invariant.
\smallskip
\begin{rem}\label{p1p2}
We recall that the smallest number in the Fukui invariant $A(f)$
is the multiplicity of $f$. Let
$$
f(x_1, \cdots , x_d) = \pm x_1^{p_1} \pm x_2^{p_2} \pm 
\cdots \pm x_d^{p_d}, \ \ 2 \le p_1 \le p_2 \le \cdots \le p_d.
$$
Then $p_1$ is determined as the smallest number in $A(f)$.
Let $n$ be the smallest number in $A(f)$
that is not divisible by $p_1$.
Suppose that $p_1$ is odd.
If $kp_1 < p_2 < (k + 1)p_1$ for some positive integer $k$,
then $n = p_2$. 
In case where $p_2 = kp_1$ for some $k$, using the argument 
of example \ref{example5} we see that $n = p_2 + 1$.
Therefore, if $kp_1 + 1 < n < (k + 1)p_1$
then $p_2 = n$.
If $n = kp_1 + 1$ then $p_2 = n - 1$ or $n$.
This implies that if $kp_1 + 1 < n < (k + 1)p_1$
then $p_2$ is determined by $A(f)$.
\end{rem}

\medskip
\begin{thm}\label{classification3} \ \ \
Let $f_i(x,y,z) = \pm x^{p_i}  \pm y^{q_i} \pm z^{r_i}$, $2 \le p_i\le q_i 
\le r_i$, $ i=1,2$, be two Brieskorn polynomials with the same Fukui 
invariants and the same zeta functions.  Then $p_1=p_2$ and one of the two 
following cases holds: 
\begin{enumerate}
\item [(i)]
$p=p_1=p_2$ is even and $f_1$ and $f_2$ belong to one of the 
following families: $\{x^p+y^{kp} +z^{kp}; k\in \N \}$, $
\{-(x^p+y^{kp} +z^{kp}); k\in \N \}$
\item [(ii)] 
$q_1=q_2$, $r_1=r_2$, and $f_1$ and $f_2$ are blow-analytically           
equivalent.
\end{enumerate}
\end{thm}

We make the 
following convention. Whenever a Brieskorn polynomial  $f(x,y,z) = \pm x^{p} 
\pm  y^{q} \pm z^{r}$ contains two terms with the same exponents and 
different signs then ``+" preceeds ``-", for instance we write 
$x^p-y^p$ instead of $-x^p +y^p$.  

\begin{rem}\label{signs}
Suppose that the $f_i$'s are written down according to the above convention.  
Then, in the second case of theorem \ref{classification3} the signs 
corresponding to the even exponents have to be the same for $i=1,2$, except 
for the case when an even exponent ($q$ or $r$) is a multiple of 
another exponent ($p$ or 
$q$) that is odd.  In the latter case the sign cannot be determined.  For 
instance we 
cannot distinguish $x^p + y^{kp}$ from  $x^p - y^{kp}$, $p$ odd, $k$ even, 
cf. proof of theorem \ref{classification}.   
\end{rem}

\begin{proof}[Proof of theorem \ref{classification3}]

Let $f(x,y,z) = \pm x^{p} \pm y^{q} \pm z^{r}$, $1<p\le q\le r$. 
We show that except the cases considered in (i) the exponents $p$, $q$, 
and $r$ are determined by the zeta functions and the Fukui invariants of 
$f$.  We  suppose that the signs in $f$ satisfy the above convention.  

First note that 
$p$ is determined by the Fukui invariant.  

If $p$ is even then the Fukui invariants with sign 
determine the sign at $x^p$ (if $p=q$ by the sign convention).   
Then, by 
corollary \ref{xmeven}, the zeta functions of $g(y,z) = \pm y^{q}
\pm z^{r}$ are  determined by the zeta functions  of $f$.  If 
$\zetaa_{g,\pm}$ are not identically equal to zero then we may use 
subsection \ref{distinction} to determine the exponents and the 
blow-analytic type of $g$.  The signs are determined as in remark 
\ref{signs}.  If the zeta functions of $g$ are identically equal 
to zero then $g(y,z) = \pm (y^{q} + z^{q})$, $q$ even.  Note that the 
Fukui invariants of $\pm x^p \pm (y^{q} + z^{q})$, $p\le q$ both even, 
are the same as the Fukui invariants of $\pm x^p \pm y^{q}$.  
The latter are given in subsection 
\ref{List}.  Thus the Fukui invariants determine $q$, and the sign, in 
all cases except (i) of the theorem.   

Suppose that $p$ is odd.  Consider the Fukui invariant $A(f)$. Let 
$n$ be the smallest number in $A(f)$ that is not divisible by $p$.  If 
$ kp+1 <n< (k+1)p$ then $q = n$.  If moreover such $q$ is even then 
$A_+(f), A_-(f)$ determine the sign at $y^q$.  
Note that if we determine the second exponent, for instance $q$ but the 
argument works also if it is $r$, so that we can determine uniquely 
the zeta functions of $\pm x^{p} \pm y^{q}$ then the remaining third 
exponent is unique by Lemma \ref{brieskorn}.  This ends the proof 
if  $ kp+1 <n< (k+1)p$.  

Suppose $n=kp+1$.  Then $q= kp$ or $kp +1$.  Consider first the case $k$ even.
Then $kp+1$ is odd.  Let 
$\mzeta_{f,\pm}= \sum_{n\ge 1} \tilde A^{pm}_n$.  By example
\ref{example2} and by Thom-Sebastiani Formula \eqref{tildeC} applied
twice to $f$, $kp+1$ equals $q$ or $r$ if and only if
$\tilde A^{\pm}_{kp+1}=0$.  If this is the case then we apply
Proposition \ref{brieskorn} to determine the remaining exponent (and the
sign as in Remark \ref{signs}).  If this is not the case then $q=kp$.
The zeta functions of  $\pm x^{p} \pm y^{q}$ do not 
depend on the signs, see example \ref{example3}, and we may apply again 
Proposition \ref{brieskorn} to determine $r$.  

Thus the only remaining case is $p$ odd, $q=kp$ or $kp+1$,
with $k$ odd.  In this case, $q$ can be determined by
the coefficients $\tilde A^+_{kp+1}$,
$\tilde A^-_{kp+1}$, $\tilde A^+_{kp+2}=\tilde A^-_{kp+2}$ of the
modified zeta function of $f$ and the Fukui invariants, that is
the knowledge whether $kp+1\in A_+(f)$ or $kp+1\in A_-(f)$.  
The computation is summarized in the table below.  
\end{proof}

\begin{table}[hbt]
\begin{tabular}{|l|c|c|c|c|c|}
\hline 
$g(y,z)=\pm y^{q} \pm z^{r}$ & $\tilde A^+_{kp+1}$ & $\tilde A^-_{kp+1}$ 
& $\tilde A^\pm _{kp+2}$ & $kp+1 \in A_+(f)$ & $kp+1 \in A_-(f)$  \\ \hline 
\hline 
$\pm y ^{kp} \pm z^{kp}$    \hfil   & -1                  &   -1 
& -1                     & yes                & yes               \\ \hline
$\pm y ^{kp} + z^{kp+1}$     & 1                  &   -1 
& -1                     & yes                & yes               \\ \hline
$\pm y ^{kp} - z^{kp+1}$       & -1                  &   1 
& -1                     & yes                & yes               \\ \hline
$\pm y ^{kp} \pm z^{kp+2}$       & 1                  &   1 
& 0                    & yes                & yes               \\ \hline
$\pm y ^{kp} \pm z^{r},r>kp+2$ & 1                  &   1 
& 1                     & yes                & yes               \\ \hline
\hline 
$y ^{kp+1} + z^{kp+1}$       & -1                  &   -1 
& -1                     & yes                & no               \\ \hline
$y ^{kp+1} - z^{kp+1}$       & 1                  &   1 
& -1                     & yes                & yes               \\ \hline
$-y ^{kp+1} - z^{kp+1}$       & -1                  &   -1 
& -1                     & no                & yes               \\ \hline
$y ^{kp+1} \pm z^{kp+2}$       & -1                  &   1 
& 0                     & yes                & no               \\ \hline
$-y ^{kp+1} \pm z^{kp+2}$       & 1                  &   -1 
& 0                     & no               & yes               \\ \hline
$y ^{kp+1} \pm z^{r},r>kp+2$       & -1                  &   1 
& 1                     & yes                & no               \\ \hline
$-y ^{kp+1} \pm z^{r},r>kp+2$       & 1                  &   -1 
& 1                     & no                & yes               \\ \hline
\end{tabular}
%\caption{}
\end{table}

\subsection{Example on blow-analytic sufficiency of jets}
\label{nB}
The zeta functions can be used to distinguish the blow-analytic types of 
functions that are not necessarily Brieskorn polynomials.  For such a 
function it may be simpler to use the standard zeta functions than the 
modified ones.  To facilitate the computations we reduce the  
Thom-Sebastiani formulae of theorem \ref{Thom-Seb} modulo 2.  Taking into 
account that always $\zetap\equiv\zetam \mod 2$, i.e. the 
coefficients satisfy $a_n^+\equiv a_n^- \mod 2$, we obtain easily  :
\begin{eqnarray}\label{T-Smod2}
& & 1+c_n^+ \equiv (1+a_n^+) (1+b_n^+) \mod 2 \\
& & 1+c_n^- \equiv (1+a_n^-) (1+b_n^-) \mod 2 . \nonumber 
\end{eqnarray}
Of course these both formulae are equivalent.   

\begin{example}
Let $\K = \R$ or $\C$. 
We consider polynomial functions 
$f_{\K}$,  $g_{\K} : (\K^3,0) \to (\K,0)$ defined by
$$
f_{\K}(x,y,z) = x^3 + xy^5 + z^3, \ \
g_{\K}(x,y,z) = x^3 + y^7 + z^3.
$$
Note that they are weighted homogeneous polynomials
with isolated singularities at $0 \in \K^3$
and the Fukui invariants of $f_{\K}$ and $g_{\K}$ 
are the same and equal $A_+ = A_- = \{ 3,4,5,\ldots \} \cup \{ \infty\}$.
Let $\phi : (\K^3,0) \to (\K,0)$ be an analytic function germ
with $j^6\phi (0) = j^6f_{\K}(0)$.
In case $\K = \R$ (resp. $\K = \C$), it follows from 
theorem \ref{FuPa} (resp. \cite{damongaffney}) that
if the Taylor expansion of $\phi$ contains a term 
of the form $ay^7$, $a \neq 0$, then $\phi$ is 
blow-analytically equivalent (resp. topologically equivalent)
to $g_{\R}$ (resp. $g_{\C}$), 
otherwise $\phi$ is 
blow-analytically equivalent (resp. topologically equivalent)
to $f_{\R}$ (resp. $f_{\C}$).
Using the formula of Milnor \& Orlik (\cite{milnororlik}),
we have $\mu (f_{\C}) = 26$ and $\mu (g_{\C}) = 24$.
Thus it follows from \cite{le} or \cite{teissier} that
$f_{\C}$ and $g_{\C}$ are not topologically equivalent.
The real jet $w = f_{\R}$ was originally given 
by W. Kucharz (\cite{kucharz}) as an example
such that $w$ is $C^0$-sufficient in $C^8$ functions as a 6-jet
but not $C^0$-sufficient in $C^7$ functions as a 7-jet.
Therefore $f_{\R}$ and $g_{\R}$ are topologically equivalent
and $w$ does not satisfy the Kuiper-Kuo condition 
even as a real 7-jet. 

We show that $f_{\R}$ and $g_{\R}$
are {\em not} blow-analytically equivalent.
As a result, $w \in J^6_{\R}(2,1)$ is not blow-analytically
sufficient.

Let us first compute $\zetaa_{f_{\R},+}(T) \mod 2$.  By \eqref{xy5}, 
$$
\zetaa_{x^3+xy^5,+} (T) \equiv \frac{T^{15}}{1-T^{15}} + \frac{T^3}{1+T^3} 
\mod 2 .
$$ 
Hence, by \eqref{T-Smod2}, the coefficients $a^+_n(f)$ of 
$\zetaa_{f_{\R},+}(T) 
\mod 2$ are given by 
\begin{equation*}
a^+_n(f) \equiv 
\begin{cases} \ 1  \mod 2 \quad \text { if }
3|n   \\
\ 0 \mod 2 \hfill \quad \text {otherwise} .
\end{cases}
\end{equation*}
A similar computation of $\zetaa_{g_{\R},+} (T) \mod 2$ shows that its 
coefficients are equal to 
\begin{equation*}
a^+_n(g) \equiv 
\begin{cases} \ 1  \mod 2 \quad \text { if }
3|n \text { or } 7|n  \\
\ 0 \mod 2 \hfill \quad \text {otherwise} .
\end{cases}
\end{equation*}
Therefore, by theorem \ref{invariance},  $f_{\R}$ and $g_{\R}$
are not blow-analytically equivalent.
\end{example}

%\end{document} 
 
%%%%%%%%%%%%%%%%%%%%%%%%%%%%%%%%%%%%%%%%%%%%%%%%%%%%%%%%%%%%%%%%%%

\end{document}